\documentclass[12pt,oneside]{amsart}
\usepackage[T1]{fontenc}

\numberwithin{equation}{section}
\newtheorem{thm}[equation]{Theorem}

\newtheorem{cor}[equation]{Corollary}

\newtheorem*{bumpfunction}{Bump function lemma}

\newtheorem{lem}[equation]{Lemma}

\newtheorem{prop}[equation]{Proposition}

\newenvironment{pf}{\proof[\proofname]}{\endproof}
\newenvironment{pf*}[1]{\proof[#1]}{\endproof}

\theoremstyle{definition}

\newtheorem{exmp}[equation]{Example}

\theoremstyle{remark}

\newtheorem*{rmk}{Remark}
\newtheorem*{rmks}{Remarks}

\usepackage{amsmath,amsthm}
\usepackage{amsfonts}

\makeatother

\begin{document}
\baselineskip=18truept


\def\C {{\mathbb C}}
\def\Cn {{\mathbb C}^n}
\def\R {{\mathbb R}}
\def\Rn {{\mathbb R}^n}
\def\Z {{\mathbb Z}}
\def\N {{\mathbb N}}
\def\cal#1{{\mathcal #1}}
\def\bb#1{{\mathbb #1}}

\def\dbar {\bar \partial }
\def\dir {{\mathcal D}}
\def\lev#1{{\mathcal L}(#1)}
\def\lap {\Delta }
\def\ol {{\mathcal O}}
\def\E {{\mathcal E}}
\def\J {{\mathcal J}}
\def\U {{\mathcal U}}
\def\V {{\mathcal V}}
\def\z {\zeta }
\def\Harm {\text {Harm}\, }
\def\grad {\nabla }
\def\dexh {\{ M_k \} _{k=0}^{\infty } }
\def\sing#1{#1_{\text {sing}}}
\def\reg#1{#1_{\text {reg}}}

\def\setof#1#2{\{ \, #1 \mid #2 \, \} }
\def\Subset{\subset\subset }

\def\holecl {M\setminus \overline M_0}
\def\hole {M\setminus M_0}

\def\nd{\frac {\partial }{\partial\nu } }
\def\ndof#1{\frac {\partial#1}{\partial\nu } }

\def\pdof#1#2{\frac {\partial#1}{\partial#2}}

\def\cinf{C^{\infty }}

\def\diam{\text {diam} \, }

\def\real{\text {Re}\, }

\def\imag{\text {Im}\, }

\def\supp{\text {supp}\, }

\def\Vol{\text {\rm vol} \, }



\def\anal{analytic }
\def\analns{analytic}

\def\bdd{bounded }
\def\bddns{bounded}

\def\cpt{compact }
\def\cptns{compact}

\def\cpx{complex }
\def\cpxns{complex}

\def\cont{continuous }
\def\contns{continuous}

\def\dime{dimension }
\def\dimens{dimension }

\def\exh{exhaustion }
\def\exhns{exhaustion}

\def\fn{function }
\def\fnns{function}

\def\fns{functions }
\def\fnsns{functions}

\def\holo{holomorphic }
\def\holons{holomorphic}

\def\mero{meromorphic }
\def\merons{meromorphic}

\def\holoconvex{holomorphically convex }
\def\holoconvexns{holomorphically convex}

\def\ircomp{irreducible component }
\def\concomp{connected component }
\def\ircompns{irreducible component}
\def\concompns{connected component}
\def\ircomps{irreducible components }
\def\concomps{connected components }
\def\ircompsns{irreducible components}
\def\concompsns{connected components}

\def\irred{irreducible }
\def\irredns{irreducible}

\def\con{connected }
\def\conns{connected}

\def\comp{component }
\def\compns{component}
\def\comps{components }
\def\compsns{components}

\def\mfld{manifold }
\def\mfldns{manifold}
\def\mflds{manifolds }
\def\mfldsns{manifolds}

\def\nbd{neighborhood }
\def\nbds{neighborhoods }
\def\nbdns{neighborhood}
\def\nbdsns{neighborhoods}

\def\harm{harmonic }
\def\harmns{harmonic}
\def\plh{pluriharmonic }
\def\plhns{pluriharmonic}
\def\plsh{plurisubharmonic }
\def\plshns{plurisubharmonic}

\def\qplsh#1{$#1$-plurisubharmonic}
\def\hplsh{$(n-1)$-plurisubharmonic }
\def\hplshns{$(n-1)$-plurisubharmonic}

\def\para{parabolic }
\def\parans{parabolic}

\def\rel{relatively }
\def\relns{relatively}

\def\str{strictly }
\def\strns{strictly}

\def\strg{strongly }
\def\strgns{strongly}

\def\cvx{convex }
\def\cvxns{convex}

\def\wrt{with respect to }
\def\wrtns{with respect to}

\def\st {such that }
\def\stns {such that}

\def\hm {harmonic measure }
\def\hmns {harmonic measure}

\def\hmib {harmonic measure of the ideal boundary of }
\def\hmibns {harmonic measure of the ideal boundary of}


\def\atil{\tilde a}
\def\btil{\tilde b}
\def\ctil{\tilde c}
\def\dtil{\tilde d}
\def\etil{\tilde e}
\def\ftil{\tilde f}
\def\gtil{\tilde g}
\def\htil{\tilde h}
\def\itil{\tilde i}
\def\jtil{\tilde j}
\def\ktil{\tilde k}
\def\ltil{\tilde l}
\def\mtil{\tilde m}
\def\ntil{\tilde n}
\def\otil{\tilde o}
\def\ptil{\tilde p}
\def\qtil{\tilde q}
\def\rtil{\tilde r}
\def\stil{\tilde s}
\def\ttil{\tilde t}
\def\util{\tilde u}
\def\vtil{\tilde v}
\def\wtil{\tilde w}
\def\xtil{\tilde x}
\def\ytil{\tilde y}
\def\ztil{\tilde z}

\def\Atil{\tilde A}
\def\Btil{\widetilde B}
\def\Ctil{\widetilde C}
\def\Dtil{\widetilde D}
\def\Etil{\widetilde E}
\def\Ftil{\widetilde F}
\def\Gtil{\widetilde G}
\def\Htil{\widetilde H}
\def\Itil{\tilde I}
\def\Jtil{\widetilde J}
\def\Ktil{\widetilde K}
\def\Ltil{\widetilde L}
\def\Mtil{\widetilde M}
\def\Ntil{\widetilde N}
\def\Otil{\widetilde O}
\def\Ptil{\widetilde P}
\def\Qtil{\widetilde Q}
\def\Rtil{\widetilde R}
\def\Stil{\widetilde S}
\def\Ttil{\widetilde T}
\def\Util{\widetilde U}
\def\Vtil{\widetilde V}
\def\Wtil{\widetilde W}
\def\Xtil{\widetilde X}
\def\Ytil{\widetilde Y}
\def\Ztil{\widetilde Z}

\def\alphatil {\tilde \alpha  }
\def\betatil {\tilde \beta  }
\def\gammatil {\tilde \gamma  }
\def\deltatil {\tilde \delta }
\def\epsilontil {\tilde \epsilon  }
\def\varepsilontil {\tilde \varepsilon  }
\def\zetatil {\tilde \zeta  }
\def\etatil {\tilde \eta  }
\def\thetatil {\tilde \theta  }
\def\varthetatil {\tilde \vartheta  }
\def\iotatil {\tilde \iota  }
\def\kappatil {\tilde \kappa  }
\def\lambdatil {\tilde \lambda  }
\def\mutil {\tilde \mu  }
\def\nutil {\tilde \nu  }
\def\xitil {\tilde \xi  }
\def\pitil {\tilde \pi  }
\def\varpitil {\tilde \varpi  }
\def\rhotil {\tilde \rho  }
\def\varrhotil {\tilde \varrho  }
\def\sigmatil {\tilde \sigma  }
\def\varsigmatil {\tilde \varsigma  }
\def\tautil {\tilde \tau  }
\def\upsilontil {\tilde \upsilon  }
\def\phitil {\tilde \phi  }
\def\varphitil {\tilde \varphi  }
\def\chitil {\tilde \chi  }
\def\psitil {\tilde \psi  }
\def\omegatil {\tilde \omega  }

\def\Gammatil {\widetilde \Gamma  }
\def\Deltatil {\tilde \Delta }
\def\Thetatil {\widetilde \Theta  }
\def\Lambdatil {\tilde \Lambda  }
\def\Xitil {\widetilde \Xi  }
\def\Pitil {\widetilde \Pi  }
\def\Sigmatil {\widetilde \Sigma  }
\def\Upsilontil {\widetilde \Upsilon  }
\def\Phitil {\tilde \Phi  }
\def\Psitil {\widetilde \Psi  }
\def\Omegatil {\widetilde \Omega  }

\def\varGammatil {\widetilde \varGamma  }
\def\varDeltatil {\tilde \varDelta  }
\def\varThetatil {\widetilde \varTheta  }
\def\varLambdatil {\tilde \varLambda  }
\def\varXitil {\widetilde \varXi  }
\def\varPitil {\widetilde \varPi  }
\def\varSigmatil {\widetilde \varSigma  }
\def\varUpsilontil {\widetilde \varUpsilon  }
\def\varPhitil {\tilde \varPhi  }
\def\varPsitil {\widetilde \varPsi  }
\def\varOmegatil {\widetilde \varOmega  }

\def\boldGammatil {\widetilde \boldGamma  }
\def\boldDeltatil {\tilde \boldDelta  }
\def\boldThetatil {\widetilde \boldTheta  }
\def\boldLambdatil {\tilde \boldLambda  }
\def\boldXitil {\widetilde \boldXi  }
\def\boldPitil {\widetilde \boldPi  }
\def\boldSigmatil {\widetilde \boldSigma  }
\def\boldUpsilontil {\widetilde \boldUpsilon  }
\def\boldPhitil {\tilde \boldPhi  }
\def\boldPsitil {\widetilde \boldPsi  }
\def\boldOmegatil {\widetilde \boldOmega  }


\def\ahat{\hat a}
\def\bhat{\hat b}
\def\chat{\hat c}
\def\dhat{\hat d}
\def\ehat{\hat e}
\def\fhat{\hat f}
\def\ghat{\hat g}
\def\hhat{\hat h}
\def\ihat{\hat i}
\def\jhat{\hat j}
\def\khat{\hat k}
\def\lhat{\hat l}
\def\mhat{\hat m}
\def\nhat{\hat n}
\def\ohat{\hat o}
\def\phat{\hat p}
\def\qhat{\hat q}
\def\rhat{\hat r}
\def\shat{\hat s}
\def\that{\hat t}
\def\uhat{\hat u}
\def\vhat{\hat v}
\def\what{\hat w}
\def\xhat{\hat x}
\def\yhat{\hat y}
\def\zhat{\hat z}

\def\Ahat{\hat A}
\def\Bhat{\widehat B}
\def\Chat{\widehat C}
\def\Dhat{\widehat D}
\def\Ehat{\widehat E}
\def\Fhat{\widehat F}
\def\Ghat{\widehat G}
\def\Hhat{\widehat H}
\def\Ihat{\hat I}
\def\Jhat{\widehat J}
\def\Khat{\widehat K}
\def\Lhat{\widehat L}
\def\Mhat{\widehat M}
\def\Nhat{\widehat N}
\def\Ohat{\widehat O}
\def\Phat{\widehat P}
\def\Qhat{\widehat Q}
\def\Rhat{\widehat R}
\def\Shat{\widehat S}
\def\That{\widehat T}
\def\Uhat{\widehat U}
\def\Vhat{\widehat V}
\def\What{\widehat W}
\def\Xhat{\widehat X}
\def\Yhat{\widehat Y}
\def\Zhat{\widehat Z}

\def\alphahat {\hat \alpha  }
\def\betahat {\hat \beta  }
\def\gammahat {\hat \gamma  }
\def\deltahat {\hat \delta }
\def\epsilonhat {\hat \epsilon  }
\def\varepsilonhat {\hat \varepsilon  }
\def\zetahat {\hat \zeta  }
\def\etahat {\hat \eta  }
\def\thetahat {\hat \theta  }
\def\varthetahat {\hat \vartheta  }
\def\iotahat {\hat \iota  }
\def\kappahat {\hat \kappa  }
\def\lambdahat {\hat \lambda  }
\def\muhat {\hat \mu  }
\def\nuhat {\hat \nu  }
\def\xihat {\hat \xi  }
\def\pihat {\hat \pi  }
\def\varpihat {\hat \varpi  }
\def\rhohat {\hat \rho  }
\def\varrhohat {\hat \varrho  }
\def\sigmahat {\hat \sigma  }
\def\varsigmahat {\hat \varsigma  }
\def\tauhat {\hat \tau  }
\def\upsilonhat {\hat \upsilon  }
\def\phihat {\hat \phi  }
\def\varphihat {\hat \varphi  }
\def\vphihat {\hat \varphi  }
\def\chihat {\hat \chi  }
\def\psihat {\hat \psi  }
\def\omegahat {\hat \omega  }

\def\Gammahat {\widehat \Gamma  }
\def\Deltahat {\hat \Delta }
\def\Thetahat {\widehat \Theta  }
\def\Lambdahat {\hat \Lambda  }
\def\Xihat {\widehat \Xi  }
\def\Pihat {\widehat \Pi  }
\def\Sigmahat {\widehat \Sigma  }
\def\Upsilonhat {\widehat \Upsilon  }
\def\Phihat {\hat \Phi  }
\def\Psihat {\widehat \Psi  }
\def\Omegahat {\widehat \Omega  }

\def\varGammahat {\widehat \varGamma  }
\def\varDeltahat {\hat \varDelta  }
\def\varThetahat {\widehat \varTheta  }
\def\varLambdahat {\hat \varLambda  }
\def\varXihat {\widehat \varXi  }
\def\varPihat {\widehat \varPi  }
\def\varSigmahat {\widehat \varSigma  }
\def\varUpsilonhat {\widehat \varUpsilon  }
\def\varPhihat {\hat \varPhi  }
\def\varPsihat {\widehat \varPsi  }
\def\varOmegahat {\widehat \varOmega  }

\def\boldGammahat {\widehat \boldGamma  }
\def\boldDeltahat {\hat \boldDelta  }
\def\boldThetahat {\widehat \boldTheta  }
\def\boldLambdahat {\hat \boldLambda  }
\def\boldXihat {\widehat \boldXi  }
\def\boldPihat {\widehat \boldPi  }
\def\boldSigmahat {\widehat \boldSigma  }
\def\boldUpsilonhat {\widehat \boldUpsilon  }
\def\boldPhihat {\hat \boldPhi  }
\def\boldPsihat {\widehat \boldPsi  }
\def\boldOmegahat {\widehat \boldOmega  }

\def\seq#1#2{\{#1_{#2}\} }

\def\aseq{ \{ a_{\nu } \} }
\def\bseq{ \{ b_{\nu } \} }
\def\cseq{ \{ c_{\nu } \} }
\def\dseq{ \{ d_{\nu } \} }
\def\eseq{ \{ e_{\nu } \} }
\def\fseq{ \{ f_{\nu } \} }
\def\gseq{ \{ g_{\nu } \} }
\def\hseq{ \{ h_{\nu } \} }
\def\iseq{ \{ i_{\nu } \} }
\def\jseq{ \{ j_{\nu } \} }
\def\kseq{ \{ k_{\nu } \} }
\def\lseq{ \{ l_{\nu } \} }
\def\mseq{ \{ m_{\nu } \} }
\def\nseq{ \{ n_{\nu } \} }
\def\oseq{ \{ o_{\nu } \} }
\def\pseq{ \{ p_{\nu } \} }
\def\qseq{ \{ q_{\nu } \} }
\def\rseq{ \{ r_{\nu } \} }
\def\sseq{ \{ s_{\nu } \} }
\def\tseq{ \{ t_{\nu } \} }
\def\useq{ \{ u_{\nu } \} }
\def\vseq{ \{ v_{\nu } \} }
\def\wseq{ \{ w_{\nu } \} }
\def\xseq{ \{ x_{\nu } \} }
\def\yseq{ \{ y_{\nu } \} }
\def\zseq{ \{ z_{\nu } \} }

\def\Aseq{ \{ A_{\nu } \} }
\def\Bseq{ \{ B_{\nu } \} }
\def\Cseq{ \{ C_{\nu } \} }
\def\Dseq{ \{ D_{\nu } \} }
\def\Eseq{ \{ E_{\nu } \} }
\def\Fseq{ \{ F_{\nu } \} }
\def\Gseq{ \{ G_{\nu } \} }
\def\Hseq{ \{ H_{\nu } \} }
\def\Iseq{ \{ I_{\nu } \} }
\def\Jseq{ \{ J_{\nu } \} }
\def\Kseq{ \{ K_{\nu } \} }
\def\Lseq{ \{ L_{\nu } \} }
\def\Mseq{ \{ M_{\nu } \} }
\def\Nseq{ \{ N_{\nu } \} }
\def\Oseq{ \{ O_{\nu } \} }
\def\Pseq{ \{ P_{\nu } \} }
\def\Qseq{ \{ Q_{\nu } \} }
\def\Rseq{ \{ R_{\nu } \} }
\def\Sseq{ \{ S_{\nu } \} }
\def\Tseq{ \{ T_{\nu } \} }
\def\Useq{ \{ U_{\nu } \} }
\def\Vseq{ \{ V_{\nu } \} }
\def\Wseq{ \{ W_{\nu } \} }
\def\Xseq{ \{ X_{\nu } \} }
\def\Yseq{ \{ Y_{\nu } \} }
\def\Zseq{ \{ Z_{\nu } \} }

\def\alphaseq { \{ \alpha _{\nu } \} }
\def\betaseq { \{ \beta _{\nu } \} }
\def\gammaseq { \{ \gamma _{\nu } \} }
\def\deltaseq { \{ \delta _{\nu } \} }
\def\epsilonseq { \{ \epsilon _{\nu } \} }
\def\varepsilonseq { \{ \varepsilon _{\nu } \} }
\def\zetaseq { \{ \zeta _{\nu } \} }
\def\etaseq { \{ \eta _{\nu } \} }
\def\thetaseq { \{ \theta _{\nu } \} }
\def\varthetaseq { \{ \vartheta _{\nu } \} }
\def\iotaseq { \{ \iota _{\nu } \} }
\def\kappaseq { \{ \kappa _{\nu } \} }
\def\lambdaseq { \{ \lambda _{\nu } \} }
\def\museq { \{ \mu _{\nu } \} }
\def\nuseq { \{ \nu  _{\nu } \} }
\def\xiseq { \{ \xi _{\nu } \} }
\def\piseq { \{ \pi _{\nu } \} }
\def\varpiseq { \{ \varpi _{\nu } \} }
\def\rhoseq { \{ \rho _{\nu } \} }
\def\varrhoseq { \{ \varrho _{\nu } \} }
\def\sigmaseq { \{ \sigma _{\nu } \} }
\def\varsigmaseq { \{ \varsigma _{\nu } \} }
\def\tauseq { \{ \tau _{\nu } \} }
\def\upsilonseq { \{ \upsilon _{\nu } \} }
\def\phiseq { \{ \phi _{\nu } \} }
\def\varphiseq { \{ \varphi _{\nu } \} }
\def\chiseq { \{ \chi _{\nu } \} }
\def\psiseq { \{ \psi _{\nu } \} }
\def\omegaseq { \{ \omega _{\nu } \} }

\def\Gammaseq { \{ \Gamma _{\nu } \} }
\def\Deltaseq { \{ \Delta _{\nu } \} }
\def\Thetaseq { \{ \Theta _{\nu } \} }
\def\Lambdaseq { \{ \Lambda _{\nu } \} }
\def\Xiseq { \{ \Xi _{\nu } \} }
\def\Piseq { \{ \Pi _{\nu } \} }
\def\Sigmaseq { \{ \Sigma _{\nu } \} }
\def\Upsilonseq { \{ \Upsilon _{\nu } \} }
\def\Phiseq { \{ \Phi _{\nu } \} }
\def\Psiseq { \{ \Psi _{\nu } \} }
\def\Omegaseq { \{ \Omega _{\nu } \} }

\def\varGammaseq { \{ \varGamma _{\nu } \} }
\def\varDeltaseq { \{ \varDelta _{\nu } \} }
\def\varThetaseq { \{ \varTheta _{\nu } \} }
\def\varLambdaseq { \{ \varLambda _{\nu } \} }
\def\varXiseq { \{ \varXi _{\nu } \} }
\def\varPiseq { \{ \varPi _{\nu } \} }
\def\varSigmaseq { \{ \varSigma _{\nu } \} }
\def\varUpsilonseq { \{ \varUpsilon _{\nu } \} }
\def\varPhiseq { \{ \varPhi _{\nu } \} }
\def\varPsiseq { \{ \varPsi _{\nu } \} }
\def\varOmegaseq { \{ \varOmega _{\nu } \} }

\def\boldGammaseq { \{ \boldGamma _{\nu } \} }
\def\boldDeltaseq { \{ \boldDelta _{\nu } \} }
\def\boldThetaseq { \{ \boldTheta _{\nu } \} }
\def\boldLambdaseq { \{ \boldLambda _{\nu } \} }
\def\boldXiseq { \{ \boldXi _{\nu } \} }
\def\boldPiseq { \{ \boldPi _{\nu } \} }
\def\boldSigmaseq { \{ \boldSigma _{\nu } \} }
\def\boldUpsilonseq { \{ \boldUpsilon _{\nu } \} }
\def\boldPhiseq { \{ \boldPhi _{\nu } \} }
\def\boldPsiseq { \{ \boldPsi _{\nu } \} }
\def\boldOmegaseq { \{ \boldOmega _{\nu } \} }

\def\amu{   a_{\mu }  }
\def\bmu{   b_{\mu }  }
\def\cmu{   c_{\mu }  }
\def\dmu{   d_{\mu }  }
\def\emu{   e_{\mu }  }
\def\fmu{   f_{\mu }  }
\def\gmu{   g_{\mu }  }
\def\hmu{   h_{\mu }  }
\def\imu{   i_{\mu }  }
\def\jmu{   j_{\mu }  }
\def\kmu{   k_{\mu }  }
\def\lmu{   l_{\mu }  }
\def\mmu{   m_{\mu }  }
\def\nmu{   n_{\mu }  }
\def\omu{   o_{\mu }  }
\def\pmu{   p_{\mu }  }
\def\qmu{   q_{\mu }  }
\def\rmu{   r_{\mu }  }
\def\smu{   s_{\mu }  }
\def\tmu{   t_{\mu }  }
\def\umu{   u_{\mu }  }
\def\vmu{   v_{\mu }  }
\def\wmu{   w_{\mu }  }
\def\xmu{   x_{\mu }  }
\def\ymu{   y_{\mu }  }
\def\zmu{   z_{\mu }  }

\def\Amu{   A_{\mu }  }
\def\Bmu{   B_{\mu }  }
\def\Cmu{   C_{\mu }  }
\def\Dmu{   D_{\mu }  }
\def\Emu{   E_{\mu }  }
\def\Fmu{   F_{\mu }  }
\def\Gmu{   G_{\mu }  }
\def\Hmu{   H_{\mu }  }
\def\Imu{   I_{\mu }  }
\def\Jmu{   J_{\mu }  }
\def\Kmu{   K_{\mu }  }
\def\Lmu{   L_{\mu }  }
\def\Mmu{   M_{\mu }  }
\def\Nmu{   N_{\mu }  }
\def\Omu{   O_{\mu }  }
\def\Pmu{   P_{\mu }  }
\def\Qmu{   Q_{\mu }  }
\def\Rmu{   R_{\mu }  }
\def\Smu{   S_{\mu }  }
\def\Tmu{   T_{\mu }  }
\def\Umu{   U_{\mu }  }
\def\Vmu{   V_{\mu }  }
\def\Wmu{   W_{\mu }  }
\def\Xmu{   X_{\mu }  }
\def\Ymu{   Y_{\mu }  }
\def\Zmu{   Z_{\mu }  }

\def\alphamu{\alpha _{\mu }}
\def\betamu{\beta _{\mu }}
\def\gammamu{\gamma _{\mu }}
\def\deltamu{\delta _{\mu }}
\def\epsilonmu{\epsilon _{\mu }}
\def\varepsilonmu{\varepsilon _{\mu }}
\def\zetamu{\zeta _{\mu }}
\def\etamu{\eta _{\mu }}
\def\thetamu{\theta _{\mu }}
\def\varthetamu{\vartheta _{\mu }}
\def\iotamu{\iota _{\mu }}
\def\kappamu{\kappa _{\mu }}
\def\lambdamu{\lambda _{\mu }}
\def\mumu{\mu _{\mu }}
\def\numu{\nu _{\mu }}
\def\ximu{\xi _{\mu }}
\def\pimu{\pi _{\mu }}
\def\varpimu{\varpi _{\mu }}
\def\rhomu{\rho _{\mu }}
\def\varrhomu{\varrho _{\mu }}
\def\sigmamu{\sigma _{\mu }}
\def\varsigmamu{\varsigma _{\mu }}
\def\taumu{\tau _{\mu }}
\def\upsilonmu{\upsilon _{\mu }}
\def\phimu{\phi _{\mu }}
\def\varphimu{\varphi _{\mu }}
\def\chimu{\chi _{\mu }}
\def\psimu{\psi _{\mu }}
\def\omegamu{\omega _{\mu }}

\def\Gammamu{\Gamma _{\mu }}
\def\Deltamu{\Delta _{\mu }}
\def\Thetamu{\Theta _{\mu }}
\def\Lambdamu{\Lambda _{\mu }}
\def\Ximu{\Xi _{\mu }}
\def\Pimu{\Pi _{\mu }}
\def\Sigmamu{\Sigma _{\mu }}
\def\Upsilonmu{\Upsilon _{\mu }}
\def\Phimu{\Phi _{\mu }}
\def\Psimu{\Psi _{\mu }}
\def\Omegamu{\Omega _{\mu }}

\def\varGammamu{\varGamma _{\mu }}
\def\varDeltamu{\varDelta _{\mu }}
\def\varThetamu{\varTheta _{\mu }}
\def\varLambdamu{\varLambda _{\mu }}
\def\varXimu{\varXi _{\mu }}
\def\varPimu{\varPi _{\mu }}
\def\varSigmamu{\varSigma _{\mu }}
\def\varUpsilonmu{\varUpsilon _{\mu }}
\def\varPhimu{\varPhi _{\mu }}
\def\varPsimu{\varPsi _{\mu }}
\def\varOmegamu{\varOmega _{\mu }}

\def\boldGammamu{\boldGamma _{\mu }}
\def\boldDeltamu{\boldDelta _{\mu }}
\def\boldThetamu{\boldTheta _{\mu }}
\def\boldLambdamu{\boldLambda _{\mu }}
\def\boldXimu{\boldXi _{\mu }}
\def\boldPimu{\boldPi _{\mu }}
\def\boldSigmamu{\boldSigma _{\mu }}
\def\boldUpsilonmu{\boldUpsilon _{\mu }}
\def\boldPhimu{\boldPhi _{\mu }}
\def\boldPsimu{\boldPsi _{\mu }}
\def\boldOmegamu{\boldOmega _{\mu }}


\def\asmu{   a^{(\mu )}  }
\def\bsmu{   b^{(\mu )}  }
\def\csmu{   c^{(\mu )}  }
\def\dsmu{   d^{(\mu )}  }
\def\esmu{   e^{(\mu )}  }
\def\fsmu{   f^{(\mu )}  }
\def\gsmu{   g^{(\mu )}  }
\def\hsmu{   h^{(\mu )}  }
\def\ismu{   i^{(\mu )}  }
\def\jsmu{   j^{(\mu )}  }
\def\ksmu{   k^{(\mu )}  }
\def\lsmu{   l^{(\mu )}  }
\def\msmu{   m^{(\mu )}  }
\def\nsmu{   n^{(\mu )}  }
\def\osmu{   o^{(\mu )}  }
\def\psmu{   p^{(\mu )}  }
\def\qsmu{   q^{(\mu )}  }
\def\rsmu{   r^{(\mu )}  }
\def\ssmu{   s^{(\mu )}  }
\def\tsmu{   t^{(\mu )}  }
\def\usmu{   u^{(\mu )}  }
\def\vsmu{   v^{(\mu )}  }
\def\wsmu{   w^{(\mu )}  }
\def\xsmu{   x^{(\mu )}  }
\def\ysmu{   y^{(\mu )}  }
\def\zsmu{   z^{(\mu )}  }

\def\Asmu{   A^{(\mu )}  }
\def\Bsmu{   B^{(\mu )}  }
\def\Csmu{   C^{(\mu )}  }
\def\Dsmu{   D^{(\mu )}  }
\def\Esmu{   E^{(\mu )}  }
\def\Fsmu{   F^{(\mu )}  }
\def\Gsmu{   G^{(\mu )}  }
\def\Hsmu{   H^{(\mu )}  }
\def\Ismu{   I^{(\mu )}  }
\def\Jsmu{   J^{(\mu )}  }
\def\Ksmu{   K^{(\mu )}  }
\def\Lsmu{   L^{(\mu )}  }
\def\Msmu{   M^{(\mu )}  }
\def\Nsmu{   N^{(\mu )}  }
\def\Osmu{   O^{(\mu )}  }
\def\Psmu{   P^{(\mu )}  }
\def\Qsmu{   Q^{(\mu )}  }
\def\Rsmu{   R^{(\mu )}  }
\def\Ssmu{   S^{(\mu )}  }
\def\Tsmu{   T^{(\mu )}  }
\def\Usmu{   U^{(\mu )}  }
\def\Vsmu{   V^{(\mu )}  }
\def\Wsmu{   W^{(\mu )}  }
\def\Xsmu{   X^{(\mu )}  }
\def\Ysmu{   Y^{(\mu )}  }
\def\Zsmu{   Z^{(\mu )}  }

\def\alphasmu{\alpha ^{(\mu )}}
\def\betasmu{\beta ^{(\mu )}}
\def\gammasmu{\gamma ^{(\mu )}}
\def\deltasmu{\delta ^{(\mu )}}
\def\epsilonsmu{\epsilon ^{(\mu )}}
\def\varepsilonsmu{\varepsilon ^{(\mu )}}
\def\zetasmu{\zeta ^{(\mu )}}
\def\etasmu{\eta ^{(\mu )}}
\def\thetasmu{\theta ^{(\mu )}}
\def\varthetasmu{\vartheta ^{(\mu )}}
\def\iotasmu{\iota ^{(\mu )}}
\def\kappasmu{\kappa ^{(\mu )}}
\def\lambdasmu{\lambda ^{(\mu )}}
\def\musmu{\mu ^{(\mu )}}
\def\nusmu{\nu ^{(\mu )}}
\def\xismu{\xi ^{(\mu )}}
\def\pismu{\pi ^{(\mu )}}
\def\varpismu{\varpi ^{(\mu )}}
\def\rhosmu{\rho ^{(\mu )}}
\def\varrhosmu{\varrho ^{(\mu )}}
\def\sigmasmu{\sigma ^{(\mu )}}
\def\varsigmasmu{\varsigma ^{(\mu )}}
\def\tausmu{\tau ^{(\mu )}}
\def\upsilonsmu{\upsilon ^{(\mu )}}
\def\phismu{\phi ^{(\mu )}}
\def\varphismu{\varphi ^{(\mu )}}
\def\chismu{\chi ^{(\mu )}}
\def\psismu{\psi ^{(\mu )}}
\def\omegasmu{\omega ^{(\mu )}}

\def\Gammasmu{\Gamma ^{(\mu )}}
\def\Deltasmu{\Delta ^{(\mu )}}
\def\Thetasmu{\Theta ^{(\mu )}}
\def\Lambdasmu{\Lambda ^{(\mu )}}
\def\Xismu{\Xi ^{(\mu )}}
\def\Pismu{\Pi ^{(\mu )}}
\def\Sigmasmu{\Sigma ^{(\mu )}}
\def\Upsilonsmu{\Upsilon ^{(\mu )}}
\def\Phismu{\Phi ^{(\mu )}}
\def\Psismu{\Psi ^{(\mu )}}
\def\Omegasmu{\Omega ^{(\mu )}}

\def\varGammasmu{\varGamma ^{(\mu )}}
\def\varDeltasmu{\varDelta ^{(\mu )}}
\def\varThetasmu{\varTheta ^{(\mu )}}
\def\varLambdasmu{\varLambda ^{(\mu )}}
\def\varXismu{\varXi ^{(\mu )}}
\def\varPismu{\varPi ^{(\mu )}}
\def\varSigmasmu{\varSigma ^{(\mu )}}
\def\varUpsilonsmu{\varUpsilon ^{(\mu )}}
\def\varPhismu{\varPhi ^{(\mu )}}
\def\varPsismu{\varPsi ^{(\mu )}}
\def\varOmegasmu{\varOmega ^{(\mu )}}

\def\boldGammasmu{\boldGamma ^{(\mu )}}
\def\boldDeltasmu{\boldDelta ^{(\mu )}}
\def\boldThetasmu{\boldTheta ^{(\mu )}}
\def\boldLambdasmu{\boldLambda ^{(\mu )}}
\def\boldXismu{\boldXi ^{(\mu )}}
\def\boldPismu{\boldPi ^{(\mu )}}
\def\boldSigmasmu{\boldSigma ^{(\mu )}}
\def\boldUpsilonsmu{\boldUpsilon ^{(\mu )}}
\def\boldPhismu{\boldPhi ^{(\mu )}}
\def\boldPsismu{\boldPsi ^{(\mu )}}
\def\boldOmegasmu{\boldOmega ^{(\mu )}}

\def\anu{   a_{\nu }  }
\def\bnu{   b_{\nu }  }
\def\cnu{   c_{\nu }  }
\def\dnu{   d_{\nu }  }
\def\enu{   e_{\nu }  }
\def\fnu{   f_{\nu }  }
\def\gnu{   g_{\nu }  }
\def\hnu{   h_{\nu }  }
\def\inu{   i_{\nu }  }
\def\jnu{   j_{\nu }  }
\def\knu{   k_{\nu }  }
\def\lnu{   l_{\nu }  }
\def\mnu{   m_{\nu }  }
\def\nnu{   n_{\nu }  }
\def\onu{   o_{\nu }  }
\def\pnu{   p_{\nu }  }
\def\qnu{   q_{\nu }  }
\def\rnu{   r_{\nu }  }
\def\snu{   s_{\nu }  }
\def\tnu{   t_{\nu }  }
\def\unu{   u_{\nu }  }
\def\vnu{   v_{\nu }  }
\def\wnu{   w_{\nu }  }
\def\xnu{   x_{\nu }  }
\def\ynu{   y_{\nu }  }
\def\znu{   z_{\nu }  }

\def\Anu{   A_{\nu }  }
\def\Bnu{   B_{\nu }  }
\def\Cnu{   C_{\nu }  }
\def\Dnu{   D_{\nu }  }
\def\Enu{   E_{\nu }  }
\def\Fnu{   F_{\nu }  }
\def\Gnu{   G_{\nu }  }
\def\Hnu{   H_{\nu }  }
\def\Inu{   I_{\nu }  }
\def\Jnu{   J_{\nu }  }
\def\Knu{   K_{\nu }  }
\def\Lnu{   L_{\nu }  }
\def\Mnu{   M_{\nu }  }
\def\Nnu{   N_{\nu }  }
\def\Onu{   O_{\nu }  }
\def\Pnu{   P_{\nu }  }
\def\Qnu{   Q_{\nu }  }
\def\Rnu{   R_{\nu }  }
\def\Snu{   S_{\nu }  }
\def\Tnu{   T_{\nu }  }
\def\Unu{   U_{\nu }  }
\def\Vnu{   V_{\nu }  }
\def\Wnu{   W_{\nu }  }
\def\Xnu{   X_{\nu }  }
\def\Ynu{   Y_{\nu }  }
\def\Znu{   Z_{\nu }  }

\def\alphanu{\alpha _{\nu }}
\def\betanu{\beta _{\nu }}
\def\gammanu{\gamma _{\nu }}
\def\deltanu{\delta _{\nu }}
\def\epsilonnu{\epsilon _{\nu }}
\def\varepsilonnu{\varepsilon _{\nu }}
\def\zetanu{\zeta _{\nu }}
\def\etanu{\eta _{\nu }}
\def\thetanu{\theta _{\nu }}
\def\varthetanu{\vartheta _{\nu }}
\def\iotanu{\iota _{\nu }}
\def\kappanu{\kappa _{\nu }}
\def\lambdanu{\lambda _{\nu }}
\def\munu{\mu _{\nu }}
\def\nunu{\nu _{\nu }}
\def\xinu{\xi _{\nu }}
\def\pinu{\pi _{\nu }}
\def\varpinu{\varpi _{\nu }}
\def\rhonu{\rho _{\nu }}
\def\varrhonu{\varrho _{\nu }}
\def\sigmanu{\sigma _{\nu }}
\def\varsigmanu{\varsigma _{\nu }}
\def\taunu{\tau _{\nu }}
\def\upsilonnu{\upsilon _{\nu }}
\def\phinu{\phi _{\nu }}
\def\varphinu{\varphi _{\nu }}
\def\chinu{\chi _{\nu }}
\def\psinu{\psi _{\nu }}
\def\omeganu{\omega _{\nu }}

\def\Gammanu{\Gamma _{\nu }}
\def\Deltanu{\Delta _{\nu }}
\def\Thetanu{\Theta _{\nu }}
\def\Lambdanu{\Lambda _{\nu }}
\def\Xinu{\Xi _{\nu }}
\def\Pinu{\Pi _{\nu }}
\def\Sigmanu{\Sigma _{\nu }}
\def\Upsilonnu{\Upsilon _{\nu }}
\def\Phinu{\Phi _{\nu }}
\def\Psinu{\Psi _{\nu }}
\def\Omeganu{\Omega _{\nu }}

\def\varGammanu{\varGamma _{\nu }}
\def\varDeltanu{\varDelta _{\nu }}
\def\varThetanu{\varTheta _{\nu }}
\def\varLambdanu{\varLambda _{\nu }}
\def\varXinu{\varXi _{\nu }}
\def\varPinu{\varPi _{\nu }}
\def\varSigmanu{\varSigma _{\nu }}
\def\varUpsilonnu{\varUpsilon _{\nu }}
\def\varPhinu{\varPhi _{\nu }}
\def\varPsinu{\varPsi _{\nu }}
\def\varOmeganu{\varOmega _{\nu }}

\def\boldGammanu{\boldGamma _{\nu }}
\def\boldDeltanu{\boldDelta _{\nu }}
\def\boldThetanu{\boldTheta _{\nu }}
\def\boldLambdanu{\boldLambda _{\nu }}
\def\boldXinu{\boldXi _{\nu }}
\def\boldPinu{\boldPi _{\nu }}
\def\boldSigmanu{\boldSigma _{\nu }}
\def\boldUpsilonnu{\boldUpsilon _{\nu }}
\def\boldPhinu{\boldPhi _{\nu }}
\def\boldPsinu{\boldPsi _{\nu }}
\def\boldOmeganu{\boldOmega _{\nu }}


\def\asnu{   a^{(\nu )}  }
\def\bsnu{   b^{(\nu )}  }
\def\csnu{   c^{(\nu )}  }
\def\dsnu{   d^{(\nu )}  }
\def\esnu{   e^{(\nu )}  }
\def\fsnu{   f^{(\nu )}  }
\def\gsnu{   g^{(\nu )}  }
\def\hsnu{   h^{(\nu )}  }
\def\isnu{   i^{(\nu )}  }
\def\jsnu{   j^{(\nu )}  }
\def\ksnu{   k^{(\nu )}  }
\def\lsnu{   l^{(\nu )}  }
\def\msnu{   m^{(\nu )}  }
\def\nsnu{   n^{(\nu )}  }
\def\osnu{   o^{(\nu )}  }
\def\psnu{   p^{(\nu )}  }
\def\qsnu{   q^{(\nu )}  }
\def\rsnu{   r^{(\nu )}  }
\def\ssnu{   s^{(\nu )}  }
\def\tsnu{   t^{(\nu )}  }
\def\usnu{   u^{(\nu )}  }
\def\vsnu{   v^{(\nu )}  }
\def\wsnu{   w^{(\nu )}  }
\def\xsnu{   x^{(\nu )}  }
\def\ysnu{   y^{(\nu )}  }
\def\zsnu{   z^{(\nu )}  }

\def\Asnu{   A^{(\nu )}  }
\def\Bsnu{   B^{(\nu )}  }
\def\Csnu{   C^{(\nu )}  }
\def\Dsnu{   D^{(\nu )}  }
\def\Esnu{   E^{(\nu )}  }
\def\Fsnu{   F^{(\nu )}  }
\def\Gsnu{   G^{(\nu )}  }
\def\Hsnu{   H^{(\nu )}  }
\def\Isnu{   I^{(\nu )}  }
\def\Jsnu{   J^{(\nu )}  }
\def\Ksnu{   K^{(\nu )}  }
\def\Lsnu{   L^{(\nu )}  }
\def\Msnu{   M^{(\nu )}  }
\def\Nsnu{   N^{(\nu )}  }
\def\Osnu{   O^{(\nu )}  }
\def\Psnu{   P^{(\nu )}  }
\def\Qsnu{   Q^{(\nu )}  }
\def\Rsnu{   R^{(\nu )}  }
\def\Ssnu{   S^{(\nu )}  }
\def\Tsnu{   T^{(\nu )}  }
\def\Usnu{   U^{(\nu )}  }
\def\Vsnu{   V^{(\nu )}  }
\def\Wsnu{   W^{(\nu )}  }
\def\Xsnu{   X^{(\nu )}  }
\def\Ysnu{   Y^{(\nu )}  }
\def\Zsnu{   Z^{(\nu )}  }

\def\alphasnu{\alpha ^{(\nu )}}
\def\betasnu{\beta ^{(\nu )}}
\def\gammasnu{\gamma ^{(\nu )}}
\def\deltasnu{\delta ^{(\nu )}}
\def\epsilonsnu{\epsilon ^{(\nu )}}
\def\varepsilonsnu{\varepsilon ^{(\nu )}}
\def\zetasnu{\zeta ^{(\nu )}}
\def\etasnu{\eta ^{(\nu )}}
\def\thetasnu{\theta ^{(\nu )}}
\def\varthetasnu{\vartheta ^{(\nu )}}
\def\iotasnu{\iota ^{(\nu )}}
\def\kappasnu{\kappa ^{(\nu )}}
\def\lambdasnu{\lambda ^{(\nu )}}
\def\musnu{\mu ^{(\nu )}}
\def\nusnu{\nu ^{(\nu )}}
\def\xisnu{\xi ^{(\nu )}}
\def\pisnu{\pi ^{(\nu )}}
\def\varpisnu{\varpi ^{(\nu )}}
\def\rhosnu{\rho ^{(\nu )}}
\def\varrhosnu{\varrho ^{(\nu )}}
\def\sigmasnu{\sigma ^{(\nu )}}
\def\varsigmasnu{\varsigma ^{(\nu )}}
\def\tausnu{\tau ^{(\nu )}}
\def\upsilonsnu{\upsilon ^{(\nu )}}
\def\phisnu{\phi ^{(\nu )}}
\def\varphisnu{\varphi ^{(\nu )}}
\def\chisnu{\chi ^{(\nu )}}
\def\psisnu{\psi ^{(\nu )}}
\def\omegasnu{\omega ^{(\nu )}}

\def\Gammasnu{\Gamma ^{(\nu )}}
\def\Deltasnu{\Delta ^{(\nu )}}
\def\Thetasnu{\Theta ^{(\nu )}}
\def\Lambdasnu{\Lambda ^{(\nu )}}
\def\Xisnu{\Xi ^{(\nu )}}
\def\Pisnu{\Pi ^{(\nu )}}
\def\Sigmasnu{\Sigma ^{(\nu )}}
\def\Upsilonsnu{\Upsilon ^{(\nu )}}
\def\Phisnu{\Phi ^{(\nu )}}
\def\Psisnu{\Psi ^{(\nu )}}
\def\Omegasnu{\Omega ^{(\nu )}}

\def\varGammasnu{\varGamma ^{(\nu )}}
\def\varDeltasnu{\varDelta ^{(\nu )}}
\def\varThetasnu{\varTheta ^{(\nu )}}
\def\varLambdasnu{\varLambda ^{(\nu )}}
\def\varXisnu{\varXi ^{(\nu )}}
\def\varPisnu{\varPi ^{(\nu )}}
\def\varSigmasnu{\varSigma ^{(\nu )}}
\def\varUpsilonsnu{\varUpsilon ^{(\nu )}}
\def\varPhisnu{\varPhi ^{(\nu )}}
\def\varPsisnu{\varPsi ^{(\nu )}}
\def\varOmegasnu{\varOmega ^{(\nu )}}

\def\boldGammasnu{\boldGamma ^{(\nu )}}
\def\boldDeltasnu{\boldDelta ^{(\nu )}}
\def\boldThetasnu{\boldTheta ^{(\nu )}}
\def\boldLambdasnu{\boldLambda ^{(\nu )}}
\def\boldXisnu{\boldXi ^{(\nu )}}
\def\boldPisnu{\boldPi ^{(\nu )}}
\def\boldSigmasnu{\boldSigma ^{(\nu )}}
\def\boldUpsilonsnu{\boldUpsilon ^{(\nu )}}
\def\boldPhisnu{\boldPhi ^{(\nu )}}
\def\boldPsisnu{\boldPsi ^{(\nu )}}
\def\boldOmegasnu{\boldOmega ^{(\nu )}}


\def\vphi {\varphi }


\def\inv{   ^{-1}  }

\def\ainv{   a^{-1}  }
\def\binv{   b^{-1}  }
\def\cinv{   c^{-1}  }
\def\dinv{   d^{-1}  }
\def\einv{   e^{-1}  }
\def\finv{   f^{-1}  }
\def\ginv{   g^{-1}  }
\def\hinv{   h^{-1}  }
\def\iinv{   i^{-1}  }
\def\jinv{   j^{-1}  }
\def\kinv{   k^{-1}  }
\def\linv{   l^{-1}  }
\def\minv{   m^{-1}  }
\def\ninv{   n^{-1}  }
\def\oinv{   o^{-1}  }
\def\pinv{   p^{-1}  }
\def\qinv{   q^{-1}  }
\def\rinv{   r^{-1}  }
\def\sinv{   s^{-1}  }
\def\tinv{   t^{-1}  }
\def\uinv{   u^{-1}  }
\def\vinv{   v^{-1}  }
\def\winv{   w^{-1}  }
\def\xinv{   x^{-1}  }
\def\yinv{   y^{-1}  }
\def\zinv{   z^{-1}  }

\def\Ainv{   A^{-1}  }
\def\Binv{   B^{-1}  }
\def\Cinv{   C^{-1}  }
\def\Dinv{   D^{-1}  }
\def\Einv{   E^{-1}  }


\def\Ginv{   G^{-1}  }
\def\Hinv{   H^{-1}  }
\def\Iinv{   I^{-1}  }
\def\Jinv{   J^{-1}  }
\def\Kinv{   K^{-1}  }
\def\Linv{   L^{-1}  }
\def\Minv{   M^{-1}  }
\def\Ninv{   N^{-1}  }
\def\Oinv{   O^{-1}  }
\def\Pinv{   P^{-1}  }
\def\Qinv{   Q^{-1}  }
\def\Rinv{   R^{-1}  }
\def\Sinv{   S^{-1}  }
\def\Tinv{   T^{-1}  }
\def\Uinv{   U^{-1}  }
\def\Vinv{   V^{-1}  }
\def\Winv{   W^{-1}  }
\def\Xinv{   X^{-1}  }
\def\Yinv{   Y^{-1}  }
\def\Zinv{   Z^{-1}  }

\def\alphainv{\alpha ^{-1}}
\def\betainv{\beta ^{-1}}
\def\gammainv{\gamma ^{-1}}
\def\deltainv{\delta ^{-1}}
\def\epsiloninv{\epsilon ^{-1}}
\def\varepsiloninv{\varepsilon ^{-1}}
\def\zetainv{\zeta ^{-1}}
\def\etainv{\eta ^{-1}}
\def\thetainv{\theta ^{-1}}
\def\varthetainv{\vartheta ^{-1}}
\def\iotainv{\iota ^{-1}}
\def\kappainv{\kappa ^{-1}}
\def\lambdainv{\lambda ^{-1}}
\def\muinv{\mu ^{-1}}
\def\nuinv{\nu ^{-1}}
\def\xiinv{\xi ^{-1}}
\def\piinv{\pi ^{-1}}
\def\varpiinv{\varpi ^{-1}}
\def\rhoinv{\rho ^{-1}}
\def\varrhoinv{\varrho ^{-1}}
\def\sigmainv{\sigma ^{-1}}
\def\varsigmainv{\varsigma ^{-1}}
\def\tauinv{\tau ^{-1}}
\def\upsiloninv{\upsilon ^{-1}}
\def\phiinv{\phi ^{-1}}
\def\varphiinv{\varphi ^{-1}}
\def\vphiinv{\varphi ^{-1}}
\def\chiinv{\chi ^{-1}}
\def\psiinv{\psi ^{-1}}
\def\omegainv{\omega ^{-1}}

\def\Gammainv{\Gamma ^{-1}}
\def\Deltainv{\Delta ^{-1}}
\def\Thetainv{\Theta ^{-1}}
\def\Lambdainv{\Lambda ^{-1}}
\def\Xiinv{\Xi ^{-1}}
\def\Piinv{\Pi ^{-1}}
\def\Sigmainv{\Sigma ^{-1}}
\def\Upsiloninv{\Upsilon ^{-1}}
\def\Phiinv{\Phi ^{-1}}
\def\Psiinv{\Psi ^{-1}}
\def\Omegainv{\Omega ^{-1}}

\def\varGammainv{\varGamma ^{-1}}
\def\varDeltainv{\varDelta ^{-1}}
\def\varThetainv{\varTheta ^{-1}}
\def\varLambdainv{\varLambda ^{-1}}
\def\varXiinv{\varXi ^{-1}}
\def\varPiinv{\varPi ^{-1}}
\def\varSigmainv{\varSigma ^{-1}}
\def\varUpsiloninv{\varUpsilon ^{-1}}
\def\varPhiinv{\varPhi ^{-1}}
\def\varPsiinv{\varPsi ^{-1}}
\def\varOmegainv{\varOmega ^{-1}}

\def\boldGammainv{\boldGamma ^{-1}}
\def\boldDeltainv{\boldDelta ^{-1}}
\def\boldThetainv{\boldTheta ^{-1}}
\def\boldLambdainv{\boldLambda ^{-1}}
\def\boldXiinv{\boldXi ^{-1}}
\def\boldPiinv{\boldPi ^{-1}}
\def\boldSigmainv{\boldSigma ^{-1}}
\def\boldUpsiloninv{\boldUpsilon ^{-1}}
\def\boldPhiinv{\boldPhi ^{-1}}
\def\boldPsiinv{\boldPsi ^{-1}}
\def\boldOmegainv{\boldOmega ^{-1}}
\title[Elementary construction of subharmonic exhaustion functions]
{Elementary construction of subharmonic exhaustion functions}
\author[T.~Napier]{Terrence Napier$^{*}$}
\address{Department of Mathematics\\Lehigh University\\Bethlehem, PA 18015}
\email{tjn2@lehigh.edu}
\thanks{$^{*}$Research partially
supported by NSF grants DMS9971462 and DMS0306441}
\author[M.~Ramachandran]{Mohan Ramachandran$^{**}$}
\address{Department of Mathematics\\SUNY at Buffalo\\Buffalo, NY 14260}
\email{ramac-m@newton.math.buffalo.edu}
\thanks{$^{**}$Research partially
supported by NSF grant DMS9626169}

 \keywords{Subharmonic function, potential theory, holomorphic line bundle}

\date{}

\begin{abstract}
By a theorem of Greene and Wu [GW], a non\cpt \con Riemannian
manifold admits a $\cinf $ \str sub\harm \exh \fnns . Demailly
provided an elementary proof of this fact in [D]. A further
simplification of Demailly's proof and some (mostly known)
applications are described. Applications include the fact that the
\holo line bundle associated to a nontrivial effective divisor on
a \cpt \con \cpx manifold $X$ admits a $\cinf $ Hermitian metric
with positive scalar curvature.
\end{abstract}

\maketitle

\section{Introduction} \label{introduction}

Let $(M,g)$ be a Riemannian manifold of dimension~$n$. The {\it
Laplace operator} $\lap _g$ for $g$ is given in local coordinates
$(x_1,\dots , x_n)$ by
$$
\lap _g\vphi =\frac{1}{\sqrt G}\sum _{i,j=1}^n\pdof{}{x_i}\biggl[
g^{ij} \sqrt G \pdof{\vphi }{x_j}\biggr] ,
$$
for every \fn $\vphi $ of class $C^2$; where
$$
G=\det (g_{ij})\qquad\text{and}\qquad (g^{ij})=(g_{ij})\inv .
$$
A $C^2$ real-valued \fn $\vphi $ is called {\it subharmonic} ({\it
\str subharmonic}) with respect to~$g$ if $\lap _g\vphi \geq 0$
(respectively, $\lap _g\vphi >0$).

A real-valued \fn $\rho $ on a topological space $X$ is called an
{\it \exh \fn} if
$$
\setof{x\in X}{\rho (x)<a} \Subset X \qquad \forall a\in \R .
$$

The main purpose of this paper is to describe a simple proof of
the following:

\begin{thm}[Greene and Wu {[}GW{]}] A \con non\cpt Riemannian
manifold $M$ admits a $\cinf $ \str sub\harm \exh \fnns .
\end{thm}

Greene and Wu actually produced a proper embedding by harmonic
\fns and obtained the above as a consequence. Thus their proof is
not elementary.   A related construction is that of Ohsawa [O] of
a strongly $n$-convex \exh \fn on an $n$-dimensional \cpx space
with no \cpt \ircompsns . Demailly [D] provided an elementary (and
relatively simple) proof of Theorem~0.1 (his proof is written for
the case of the Laplace operator of a Hermitian metric, but it can
be modified to give the above theorem). His method is a version of
the classical idea in Runge theory in one \cpx variable of pushing
singularities to infinity using a local construction. His local
construction, although short, requires some calculations which are
not completely transparent. Since sub\harm \fns are fundamental
objects, it is natural to search for a construction which is as
simple as possible.

In this paper we consider another proof along these lines, in
which the local construction is very simple and transparent. The
fundamental observation is that one can produce a $\cinf $ bump
\fn which is sub\harm outside an arbitrarily small set:

\begin{bumpfunction}
Let $B$ be a domain in $M$, let $K$ be a \cpt subset of $B$, and
let $W$ be a nonempty open subset of $B\setminus K$. Then there
exists a nonnegative $\cinf $ \fn $\alpha $ on $M$ \st $\alpha
\equiv 0$ on $M\setminus B$, $\alpha >0$ and $\lap \alpha >0$ on
$K$, and $\lap \alpha \geq 0$ on $M\setminus W$.
\end{bumpfunction}
\begin{pf*}{Sketch of the proof}
We may assume without loss of generality that $W\Subset B\Subset
M$ and we may fix a domain~$U$ and a nonnegative $\cinf $ \fn
$\rho $ on $M$ \st
$$
K\cup \overline W\subset U\Subset B,
\text{ }
\rho >0 \text{ on } \overline U, \text { and }
\rho < 0 \text{ on } M\setminus B.
$$
Replacing $\rho $ by an approximating Morse \fn (see, for example,
[GG]), we may also
assume that $\rho $ has only isolated critical points in $B$. Fix
a regular value $\epsilon >0$ for $\rho $ with $\rho >\epsilon $
on $\overline U$ and let $V$ be the \comp of $\setof{x\in M}{\rho
(x)>\epsilon }$ containing $U$. Thus $U\Subset V\Subset B$.

We will say that a mapping $\Phi : N\to N$ of a \con smooth
manifold $N$ of dimension~$\geq 2$ onto itself has {\it \cpt
support} if $\Phi $ is equal to the identity outside a \cpt set.
Given two points $p,q\in N$, there exists a $\cinf $
diffeomorphism $\Phi : N\to N$ with \cpt support \st $\Phi (p)=q$
(for the set of points $q$ in $N$ to which $p$ can be moved by
such a diffeomorphism is open and closed). For distinct points
$p_1,\dots , p_m,q_1,\dots ,q_m$ in $N$, one gets such a $\Phi $
with $\Phi (p_j)=q_j$ for $j=1,\dots , m$ by forming a compactly
supported diffeomorphism $\Phi _j$ of $N\setminus \{ p_1,\dots
,\hat p_j ,\dots , p_m, q_1,\dots ,\hat q_j ,\dots ,q_m \} $
moving $p_j$ to $q_j$ for each $j=1,\dots ,m$ and letting $\Phi $
be the composition of the extensions by the identity for $\Phi
_1,\dots ,\Phi _m$.

Thus we may move into $W$ the critical points of $\rho $ in $V$ by
a diffeomorphism with \cpt support in~$V$, and hence we may assume
that $\nabla \rho \neq 0$ at each point in $\overline V\setminus W
\supset K$. For $R>0$, let $\beta \equiv e^{R\rho }-e^{R\epsilon
}$. Then
$$
\lap \beta = Re^{R\rho }(\lap \rho +R|\nabla \rho |^2)>0
$$
on $\overline V\setminus W$, provided $R\gg 0$.
Finally, fixing a $\cinf $ \fn
$\chi :\R \to \R $ \st $\chi (t)=0$ for $t\leq 0$
and
$\chi '(t)>0$ and $\chi ''(t)\geq 0$ for $t>0$, we
get a nonnegative $\cinf $ \fn
$$
\alpha \equiv \left\{
\begin{aligned}
&\chi (\beta )
& \quad \text{on } V\\
&0 &\quad \text {on } M\setminus V
\end{aligned}
\right.
$$
On $V\setminus W$, we have $\alpha >0$ and
$$
\lap \alpha =\chi '(\beta ) \lap \beta
+\chi ''(\beta )|\nabla \beta |^2>0.
$$
It follows that $\alpha $ has the required properties.
\end{pf*}

For $B$ a coordinate $n$-dimensional rectangle (or another nice
set), one can easily construct such a \fn $\alpha $ explicitly
(see Lemmas~1.6--1.8). Moreover, such bump \fns are all that are
needed to
reduce the construction of a \str sub\harm \exh \fn to point set
topology (so the proof is very elementary).

One pushes the bad set off to infinity (in the usual way) as
follows. Given a point $p\in M$, there is a locally finite
sequence of \rel \cpt domains $\seq B\nu _{\nu =1}^{\infty }$ with
$$
p\in B_1\text{ and } B_\nu \cap B_{\nu +1}\neq \emptyset \text {
for } \nu =1,2,3, \dots .
$$
Hence there exist nonempty disjoint
open sets $\seq W\nu _{\nu =0}^{\infty }$
\st $p\in N_p\equiv W_0\Subset B_1$ and
$W_\nu \Subset B_\nu \cap B_{\nu +1} \text { for }
\nu =1,2,3, \dots $ and, as in the lemma, $\cinf $ bump \fns $\seq
\alpha\nu _{\nu =1}^{\infty }$  \stns , for each $\nu =1,2,3,
\dots $, we have $\supp \alpha _\nu \subset B_\nu $, $\alpha _\nu
>0$ and $\lap\alpha _\nu >0$ on $\overline W_{\nu -1}$, and
$\lap\alpha _\nu \geq 0$ on $M\setminus W_{\nu }$. For constants
$0< r_1 \ll r_2 \ll r_3 \ll \cdots $, we get a $\cinf $
sub\harm \fn
$$
\beta _p\equiv \sum _{\nu =1}^\infty r_\nu \alpha _\nu
$$
with $\supp \beta _p\subset Q_p \equiv \bigcup _{\nu =1}^\infty
B_\nu $ and $\beta _p >0 $ and $\lap\beta _p >0$ on the \nbd $N_p$
of~$p$. Paracompactness implies that we can form a locally finite
covering $\seq N{p_j}$ of $M$ by such sets and a corresponding
locally finite collection $\seq Q{p_j}$. Thus, for $R_j\gg 0$ for
$j=1,2,3,\dots $, we get a $\cinf $ \str sub\harm \exh \fn
$$
\vphi \equiv \sum _{j=1}^\infty R_j \beta _{p_j}.
$$
In fact, the above arguments actually give the following analogue
of Urysohn's lemma:
\begin{thm}[cf.~Theorem~1.13]
Suppose $U$ is a domain in a \con non\cpt Riemannian manifold $M$,
$C$ is a \con non\cpt closed subset of $M$ with $C\subset U$, and
$\rho $ is a positive \cont \fn on $M$. Then there exists a
nonnegative $\cinf $ sub\harm \fnns~$\vphi $ on $M$ \st
$\vphi \equiv 0$ on $M\setminus U$ and $\vphi >\rho $
and $\lap \vphi >\rho $ on $C$.
\end{thm}
\begin{rmk}
The existence of such a set $C\subset U$ is a necessary condition
(see the remarks following Theorem~1.13). In the terminology of
[EM], $U$ has an {\it exit to}~$\infty $ (relative to~$M$).
\end{rmk}

A detailed proof of Theorem~0.1 (in fact, a proof of the existence
of an exhausting subsolution for a more general elliptic operator)
appears in Section~1. Some (mostly known) applications are
described in Section~2. These include the fact that the \holo line
bundle associated to a nontrivial effective divisor on a \cpt \con
\cpx manifold $X$ admits a $\cinf $ Hermitian metric with positive
scalar curvature (Theorem~2.3).

\section{Construction of exhausting subsolutions}

Throughout this section, $M$ will denote a \con smooth manifold of
dimension~$n$ and $A$ will denote a second order locally uniformly
elliptic linear differential operator with locally bounded
coefficients.
Thus, in local coordinates $(x_1,\dots ,x_n)$,
$$
A=\sum _{i,j=1}^n a_{ij}\frac{\partial ^2}{\partial x_i\partial
x_j} +\sum _{i=1}^n b_i \pdof {}{x_i} +c;
$$
where $a_{ij}$ for each $i$~and~$j$, $b_i$ for each~$i$, and~$c$
are locally bounded real-valued  \fns and $(a_{ij})$ is a
symmetric matrix-valued function whose eigenvalues are locally
bounded below by a positive constant.

Theorem~0.1 is a special case of the following:

\begin{thm}
If $M$ is non\cpt and $\rho $ is a positive \cont \fn on $M$, then
there exists a $\cinf $ \fn $\vphi $ on $M$ \st $\vphi >\rho $ and
 $A\vphi >\rho $.
\end{thm}

The main step in the proof is the following:

\begin{prop}
Suppose $K$ is a \cpt subset of $M$, $U$ is a \comp of $M\setminus
K$  which is not relatively \cpt in $M$, and $p\in U$. Then there
exists a $\cinf $ \fn $\alpha $ \st
\begin{enumerate}
\item[(i)] $\alpha \geq 0$ and $A\alpha \geq 0$ on $M$,
\item[(ii)] $\supp \alpha \subset U$,
\item[(iii)] $\alpha (p)>0$, and
\item[(iv)] $A\alpha >1$ on a \nbd of~$p$.
\end{enumerate}
\end{prop}

\begin{rmk}
Since the coefficients of the operator $A$ are only locally
bounded, the condition~(iv) is stronger than the condition
$A\alpha (p)>1$.
\end{rmk}

The following equivalent version implies that a \cpt set which is
{\it topologically Runge} is \cvx with respect to \fns $\alpha $
satisfying $A\alpha \geq 0$:

\begin{prop}
Let $K$ be a \cpt subset of $M$ whose complement has no \rel \cpt
\compsns . Then, for each point $p\in M\setminus K$, there is a
$\cinf $ nonnegative \fn $\alpha $ on $M$ \st $A\alpha \geq 0$
on~$M$, $\alpha \equiv 0$ on $K$, $\alpha (p)>0$, and $A\alpha >1$
near $p$.
\end{prop}

\begin{rmk}
If the coefficients are (for example) $C^1$ and the constant term
$c\leq 0$ (for example, if $A=\lap $), then a nonconstant
subsolution on a domain cannot attain a positive maximum and,
therefore, the converse will also hold. That is, if such a \fn
$\alpha $ exists for some point $p\in M\setminus K$, then
the \comp of $M\setminus K$ containing $p$
is not \rel \cptns .
\end{rmk}

 Proposition~1.2 and Proposition~1.3 together with standard arguments
in Runge theory give Theorem~1.1. Proofs are provided for the
convenience of the reader. For this, we need two elementary
observations (cf. Malgrange~[M] or Narasimhan~[N]).

\begin{lem}
Let $X$ be a non\cptns , \conns , locally \conns , locally \cptns , Hausdorff
topological space. If $K$ is a \cpt subset of $X$ and $\widehat K$
is the union of $K$ with all of the \rel \cpt \comps of
$X\setminus K$, then $\widehat K$ is \cptns , $X\setminus \widehat
K$ has only finitely many \compsns , and each
\comp of $X\setminus
\widehat K$ has non\cpt closure.
\end{lem}
\begin{pf}
We may assume without loss of generality that
$K\neq \emptyset $. Since $X$ is Hausdorff, $K$ is closed and,
since $X$ is locally \conns , the \comps of $X\setminus K$ are
open. It follows that $\widehat K$ is a closed set whose
complement has no \rel \cpt \comps (since $X\setminus \widehat K$
is the union of \comps of $X\setminus K$ with non\cpt closure).

Since $X$ is locally \cpt Hausdorff, we may choose a \rel \cpt
\nbd $\Omega $ of $K$ in $X$. The \comps of $X\setminus K$ are
open and disjoint, so only finitely many meet the \cpt set
$\partial \Omega \subset X\setminus K$.  By replacing $\Omega $ by
the union of $\Omega $ with all \rel \cpt \comps of $X\setminus K$
meeting $\partial \Omega $, we may assume that no \rel \cpt \comp
of $X\setminus K$ meets $\partial \Omega $. On the other hand,
every \compns~$E$ of $X\setminus K$ must satisfy
$$
\overline E\cap K=\partial E\neq \emptyset .
$$
For $E$ is open and closed relative to $X\setminus K$, so
$\partial E\subset K$, while $E\neq X$, so $\partial E=\overline
E\setminus E\neq \emptyset $ ($E$ cannot be both open and closed
in the \con space $X$). It follows that, if $E$ meets $X\setminus
\Omega $, then $E$ meets $\partial \Omega $ and hence $E$ is {\it
not} \rel \cpt in $X$. Thus
$$
X\setminus \Omega \subset E_1\cup \cdots \cup E_m
$$
for finitely many \comps $E_1,\dots ,E_m$ of $X\setminus K$, none
\rel \cpt in $X$, and $\widehat K\subset \Omega \Subset X$.  The
claim now follows.
\end{pf}
\begin{lem}
Let $X$ be a second countable, non\cptns ,
\conns , locally \conns , locally
\cptns , Hausdorff topological space. Then there is a sequence of
\cpt sets $\seq K\nu _{\nu =1}^\infty $ \st $X=\bigcup _{\nu
=1}^{\infty }K_\nu $ and, for each~$\nu $, $K_\nu \subset \overset
\circ K_{\nu +1}$ and $\widehat K_\nu =K_\nu $, where $\widehat
K_\nu $ is defined as in Lemma~1.4.
\end{lem}
\begin{pf}
We may inductively choose a sequence of \cpt sets $\seq {K'}\nu $
\st $X=\bigcup _{\nu =1}^{\infty }K'_\nu $ and, for each~$\nu $,
$K'_\nu \subset \overset {\circ \quad}{K'_{\nu +1}}$.
Setting $K_\nu =\widehat {K'_\nu }$ for each~$\nu $ yields the
desired sequence.
\end{pf}

\begin{pf*}{Proof of Theorem~1.1}
By Lemma~1.5, we may choose a sequence of nonempty \cpt sets $\seq
K\nu $ \st $M=\bigcup _{\nu =1}^{\infty }K_\nu $ and, for
each~$\nu $, $K_\nu \subset \overset \circ K_{\nu +1}$ and
$M\setminus K_\nu $ has no \rel \cpt \compsns . Set $K_0=\emptyset
$.

Given $p\in M$, there is a unique $\nu =\nu (p)$ with $p\in K_{\nu
+1}\setminus K_\nu $ and we may apply Proposition~1.3 to get a
$\cinf $ nonnegative \fn $\alpha _p$ and a \rel \cpt \nbd $V_p$
of~$p$ in $M\setminus K_\nu $ \st $A\alpha _p\geq 0$ on $M$,
$\alpha _p\equiv 0$ on $K _\nu $, and $\alpha _p>\rho $ and
$A\alpha _p>\rho $ on $V_p$ (one obtains the last two conditions
by multiplying by a sufficiently large positive constant). Thus we
may choose a sequence of points $\seq pk$ in $M$ and corresponding
\fns $\seq \alpha {p_k}$ and \nbds $\seq V{p_k}$ so that $\seq
V{p_k}$ forms a locally finite covering of~$M$ (for example, one
may take $\seq pk$ to be an enumeration of the countable set $\cup
_{\nu =0}^\infty Z_\nu $ where, for each $\nu $, $Z_\nu $ is a
finite set of points in $M\setminus \overset \circ K_\nu $ \st $\{
V_p\} _{p\in Z_\nu }$ covers $K_{\nu +1}\setminus \overset \circ
K_\nu $). The collection $\{ \supp \alpha _{p_k} \} $ is then
locally finite in $M$ since $\supp \alpha _{p_k}\subset M\setminus
K_\nu $ whenever $p_k\not\in K_\nu $. Hence the sum $\sum
_{k=1}^\infty \alpha _{p_k}$ is locally finite and, therefore,
convergent to a $\cinf $ \fn $\vphi $ on $M$ satisfying $\vphi
\geq \alpha _{p_k}
>\rho$ and  $A\vphi \geq A\alpha _{p_k}
>\rho $ on $V_{p_k}$ for each~$k$. Therefore, since $\seq V{p_k}$
covers $M$, we get $\vphi
>\rho $ and $A\vphi
>\rho $ on $M$.
\end{pf*}

It remains to prove Proposition~1.2.

\begin{lem}
Each point $p\in M$ has a \rel \cpt \con \nbdns~$V$ \stns, for
each point $q\in V$, there is a $\cinf $ nonnegative \fnns~$\rho $
on $M$ \st $\rho \equiv 0$ on $M\setminus V$, $\rho >0$ on $V$,
and $q$ is the unique critical point of $\rho $ in $V$ (hence
$\rho (q)=\max \rho $).
\end{lem}
\begin{pf}
We may assume without loss of generality that $M$ is an open
subset of $\R ^n$, $p$ is in the cube $V =(-1,1)\times \cdots
\times (-1,1)$, and $V\Subset M$. Let $q=(a_1,\dots ,a_n)\in V$
and, for each $i=1, \dots , n$, fix a $\cinf $ \fn $\lambda _i :\R
\to [0,\infty )$ \st $\lambda _i \equiv 0$ on $\R\setminus
(-1,1)$, $\lambda _i>0 $ on $(-1,1)$, and $a_i$ is the unique
critical point of $\lambda _i $ in $(-1,1)$; for example, the \fn
$$
\lambda _i(t) = \left\{
\begin{aligned}
&{\exp} \biggl( -\frac {(t-a_i)^2}{1-t^2} \biggr)
& \quad \text{if } |t|<1\\
&0 &\quad \text {if } |t|\geq 1
\end{aligned}
\right.
$$
The \fn $\rho $ given by
$$
\rho (x) =\prod _{i=1}^n \lambda _i(x_i) \qquad \forall \,
x=(x_1,\dots ,x_n)\in M
$$
then has the required properties.
\end{pf}

\begin{lem}
Let $\vphi : M\to (r,s)\subset \R $ be a $C^2$ \fnns , let $K$ be
a \cpt subset of $M$ which does not contain any critical points
for $\vphi $, and let $\chi :(r,s)\to \R $ be a $C^2$ \fn
satisfying $\chi '' \geq |\chi '|$ and $\chi ''\geq |\chi |$.
Then, for every $\epsilon $ and $R$ with $1\gg \epsilon >0$ and
$R\gg 0$, we have $A[\chi (R\vphi )]\geq \epsilon R^2 \chi
''(R\vphi )$ on $K$.
\end{lem}
\begin{pf}
Locally, we have
$$
A=\sum _{i,j=1}^n a_{ij}\frac{\partial ^2}{\partial x_i\partial
x_j} +\sum _{i=1}^n b_i \pdof {}{x_i} +c
$$
(with $a_{ij}=a_{ji}$). Hence, for $R>0$, we have
\begin{align*}
A[\chi (R\vphi)]&=\sum  a_{ij}\chi '(R\vphi )R\frac{\partial
^2\vphi } {\partial x_i\partial x_j}+\sum  a_{ij}\chi ''(R\vphi
)R^2\pdof{\vphi }{x_i}\pdof{\vphi }{x_j}\\
&\qquad\qquad\qquad\qquad\qquad\qquad\qquad\qquad +\sum  b_i \chi
'(R\vphi )R\pdof {\vphi }{x_i} +c\chi (R\vphi ) \\&= R^2\chi
''(R\vphi )\sum a_{ij}\pdof{\vphi }{x_i}\pdof{\vphi }{x_j}+ R \chi
'(R\vphi )\biggl[\sum a_{ij}\frac{\partial ^2\vphi } {\partial
x_i\partial x_j} +\sum b_i \pdof {\vphi }{x_i} \biggr] +c\chi
(R\vphi ).
\end{align*}
Since $A$ is locally uniformly elliptic with locally bounded
coefficients and $d\vphi \neq 0$ at each point in the \cpt set
$K$, it follows that there exist constants $\delta >0$ and $N>0$
(which do not depend on $R$) \stns , at each point in~$K$,
$$
A[\chi (R\vphi)]\geq R^2\delta \chi ''(R\vphi )-R|\chi '(R\vphi )|
N-N|\chi (R\vphi )| \geq \chi ''(R\vphi )(\delta R^2-RN-N).
$$
We have $\delta R^2-RN-N>\epsilon R^2$ for $\frac 12\delta >
\epsilon
>0$ and $R\gg 0$, so the claim follows.
\end{pf}

\begin{lem}
There exists a $\cinf $ \fn $\chi :\R \to \R $ \st
\begin{enumerate}
\item[(i)] $\chi (t)=0$ for $t\leq 0$, and
\item[(ii)] $\chi ''(t)\geq \chi '(t)\geq \chi (t)>0$ for $t>0$.
\end{enumerate}
\end{lem}
\begin{pf}
For example, if $a,b\geq 1$, then the $\cinf $ \fn
$$
\chi (t) = \left\{
\begin{aligned}
&{\exp} \bigl( at-(b/t) \bigr)
& \quad \text{if } t>0\\
&0 &\quad \text {if } t\leq 0
\end{aligned}
\right.
$$
satisfies (i) and (ii).
\end{pf}

Lemmas~1.6--1.8 allow one to produce bump \fns which are
subsolutions outside a small set. To push the bad set off to
infinity, we require chains of such bump \fnsns . For this, we
recall an elementary fact from point set
topology. It is convenient and
instructive to have this fact in a form which is slightly stronger
than is needed at present.

\begin{lem}
Let $X$ be a \conns , locally \conns , locally \cptns , Hausdorff
topological space, let $\mathcal B$ be a countable collection of
\con open subsets which is a basis for the topology in $X$, and
let $U$ be a \con open subset which is not \rel \cpt in $X$.
Suppose that there exists a \con non\cpt closed subset of $X$
which is contained in $U$. Then
\begin{enumerate}
\item[(i)] For any \con non\cpt closed subset $C$ of $X$ with $C\subset U$,
there exists a sequence of \con open subsets $\seq U\nu $ of $X$
\st $C\subset U_1$, $U=\bigcup _{\nu =1}^\infty U_\nu $, and, for
each $\nu $, $\overline U_\nu $ is non\cpt and $\overline U_\nu
\subset U_{\nu +1}$; and
\item[(ii)] For each point $p\in U$, there is a locally finite
sequence of basis elements $\seq Bj$ \st $p\in B_1$ and, for
each~$j$, $B_j\Subset U$ and $B_j\cap B_{j+1}\neq \emptyset $.
\end{enumerate}
If, in addition, $X$ is locally path \conns , then
\begin{enumerate}
\item [(iii)] For each point $p\in U$, there is a proper \cont map
$\gamma :[0,\infty )\to X$ with $\gamma (0)=p$ and $\gamma
([0,\infty ))\subset U$ (i.e.~a path in $U$ from $p$ to~$\infty
$).
\end{enumerate}
\end{lem}
\begin{rmk} Conversely, each of the properties
(ii)~and~(iii) clearly implies the existence of a \con non\cpt
closed subset of $X$ which is contained in $U$.
\end{rmk}
\begin{pf}
We first observe that there is a sequence of \con open subsets
$\seq \Omega\nu $ of~$U$ \st $U=\bigcup _{\nu =1}^\infty \Omega
_\nu $ and, for each~$\nu $, $\Omega _\nu \Subset \Omega _{\nu
+1}$. For we may choose a covering of $U$ by a sequence of basis
elements $\seq Gj$ which are \rel \cpt in $U$. For each $\nu $,
let $\Gamma _\nu $ be the \concomp of $G_1\cup \cdots \cup G_\nu $
containing $G_1$. Then $\Gamma \equiv \bigcup \Gamma _\nu $ is
equal to $U$. For if $p\in \overline\Gamma \cap U$, then $p\in
G_j$ for some $j$ and $G_j$ must meet $\Gamma _\mu $ for some~$\mu
$. Therefore, for $\nu
>\max (j,\mu )$, $G_j\cup \Gamma _\mu $ is a \con subset of
$G_1\cup \cdots \cup G_\nu $ containing $G_1$ and hence
$$
p\in G_j\subset \Gamma _\nu \subset \Gamma .
$$
Thus $\Gamma $ is both open and closed relative to $U$ and is
therefore equal to $U$. A suitable subsequence of $\seq \Gamma\nu
$ (chosen inductively) will then have each term \rel \cpt in the
next term, as required. We may also choose a sequence of open
subsets $\seq\Theta\nu $ \st $X=\bigcup _{\nu =1}^{\infty
}\Theta_\nu $ and, for each~$\nu $, $\Theta_\nu \Subset\Theta
_{\nu +1}$ (the above gives a proof of this elementary fact which
was used in the proof of Lemma~1.5). Set $\Omega _0=\Theta
_0=\Theta _{-1}=\emptyset $.

Next, we observe that, for any set $C$ as in~(i), there is a
countable locally finite (in $X$) covering $\mathcal A_C$ of $C$
by basis elements which are \rel \cpt in $U$. For we may take
$\mathcal A_C=\bigcup _{\nu =1}^\infty \mathcal A_C^{(\nu )}$,
where, for each $\nu =1,2,3,\dots $, $\mathcal A_C^{(\nu )}$ is a
finite covering of the \cpt set
$C\cap (\overline \Theta _{\nu
}\setminus \Theta _{\nu -1})$ by basis elements which are \rel \cpt
in $U\setminus \overline \Theta _{\nu -2}$.

For the proof of (i), we may choose the sequence
$\seq\Omega \nu $ so
that $\Omega _1\cap C\neq \emptyset $. Let $C_0=C$ and $\Omega
_0=U_0=\emptyset $. Given \con open sets $U_0,\dots , U_\nu $ and
\con closed sets $C_0,\dots ,C_\nu $ \stns , for $\mu =1,\dots
,\nu $, we have
$$
\overline \Omega _{\mu -1}\cup C_{\mu -1}\subset U_\mu \subset
\overline U_\mu =C_\mu \subset U
$$
(which holds vacuously if $\nu =0$), we may choose $U_{\nu +1}$ to
be the union of those elements of the collection $\mathcal
A_{\overline \Omega _\nu \cup  C_\nu }$ which meet the \con
non\cpt closed set  $\overline \Omega _\nu \cup  C_\nu $ and set $
C_{\nu +1}=\overline U_{\nu +1}$. Proceeding, we get a sequence
$\seq U\nu $ with the required properties.

For the proof of (ii), we may fix a \con non\cpt closed subset $C$
of $X$ with $p\in C\subset U$ (for this, we may take $\seq U\nu $
as in~(i) and let $C=\overline U_\nu $ for some $\nu \gg 0$). For
each point $q\in C$, there is a finite sequence of elements
$B_1,\dots ,B_k$ of $\mathcal A_C$ which forms a {\it chain from}
$p$ {\it to} $q$; that is, $p\in B_1$, $q\in B_k$, and $B_j\cap
B_{j+1}\neq \emptyset $ for $j=1,\dots , k-1$ (we will call $k$
the {\it length} of the chain). For the set $E$ of points $q$ in
$C$ for which there is a chain from $p$ to $q$ is clearly nonempty
and open relative to $C$. On the other hand, $E$ is also closed
because, if $q\in \overline E$, then $q\in B$ for some set $B\in
\mathcal A_C$ and there must be some point $r\in B\cap E$. A chain
$B_1,\dots , B_k$ from $p$ to~$r$ yields the chain $B_1,\dots ,
B_k,B$ from $p$ to~$q$. Thus $E=C$. Observe that if $q\in E$ and
$B_1,\dots , B_k$ is a chain of minimal length from $p$ to~$q$,
then the sets $B_1,\dots ,B_k$ are distinct.

Now since $C$ is non\cpt and closed, we may choose a sequence of
points $\seq q\nu$ in $C$ with $\qnu \to \infty $ in~$X$ (for
example, $\qnu \in C\setminus \Thetanu $ for each~$\nu $) and, for
each $\nu $, we may choose a chain $\Bsnu _1, \dots ,\Bsnu
_{k_{\nu }}$ of minimal length from $p$ to~$\qnu $. Since the
elements of $\mathcal A_C$ are \rel \cpt in $U$ and $\mathcal A_C$
is locally finite in $X$, there are only finitely many possible
choices for $\Bsnu _j$ for each~$j$ (only finitely many elements
of $\mathcal A_C$ will be in some chain of length $j$ from $p$).
Moreover, for each fixed $j\in \N $, we have $k_\nu >j$ for $\nu
\gg 0$, because the set of points in $C$ joined to $p$ by a chain
of length~$\leq j$ is \rel \cpt in $C$ while $\qnu \to \infty $.
Therefore, after applying a diagonal argument and passing to the
associated subsequence of $\seq q\nu $, we may assume that, for
each~$j$, there is an element $B_j\in \mathcal A_C$ with $\Bsnu
_j=B_j$ for all $\nu \gg 0$. Thus we get an infinite chain of
distinct elements $\seq Bj$ from $p$ to infinity as required
in~(ii) (local finiteness in $X$ is guaranteed since $\mathcal
A_C$ is locally finite and the elements $\seq Bj$ are distinct).

Finally, suppose $X$ is locally path connected. Then the sets
$\seq Bj$ as in~(ii) are path \con and, setting $p_0=p$ and,
choosing $p_j\in B_j\cap B_{j+1}$ for each $j=1,2,3,\dots $, we
may take $\gamma |_{[j-1,j]}$ to be a path in $B_j$ from $p_{j-1}$
to $p_j$ for each~$j$.
\end{pf}

\begin{pf*}{Proof of Proposition 1.2}
We first observe that, if $V$ is a set with the properties
described in Lemma~1.6, $D$ is a \cpt subset of~$V$, and $W$ is a
nonempty open subset of $V\setminus D$, then there is a
nonnegative $\cinf $ \fn $\beta $ with \cpt support in~$V$ \st
$A\beta \geq 0$ on $M\setminus W$ and $\beta >0$ and $A\beta
>1$ on $D$. For we may choose a point $q\in W$, a $\cinf $
nonnegative \fn $\rho $ on $M$ which is positive on $V$ and has
unique critical point $q$ in $V$, a $\cinf $ \fn $\chi $ on~$\R $
as in Lemma~1.8, and a constant $\epsilon >0$ with $\rho >\epsilon
$ on~$D$.  By Lemma~1.7 (applied to the \cpt set $K=\setof{x\in
M\setminus W}{\rho (x)\geq \epsilon }\subset V\setminus W$), for
$R\gg 0$, the \fn $\beta \equiv \chi (R(\rho -\epsilon ))$ will
have the required properties.

Next, by Lemma~1.9, given a point $p\in U$, there is a locally
finite (in $X$) sequence of \rel \cpt open subsets $\seq Vm$ of
$U$ \st $p\in V_1$ and, for each~$m$, $V_m$ has the properties
described in Lemma~1.6 and $V_m\cap V_{m+1}\neq \emptyset $. Hence
we may choose a sequence of disjoint nonempty open sets ${\seq
Wm}_{m=0}^\infty $ \st $p\in W_0\Subset V_1$ and, for each $m\geq
1$, $W_m\Subset V_m\cap V_{m+1}$.

By the first observation, there is a sequence of nonnegative
$\cinf $ \fns ${\seq \beta m}_{m=1}^\infty $ \stns , for each~$m$,
$\beta _m$ is compactly supported in $V_m$, $A\beta _m\geq 0$ on
$M\setminus W_m$, and $\beta _m>0$ and $A\beta _m>1$ on $\overline
W_{m-1}$.  We will choose positive constants $\seq Rm$ inductively
so that, for each $m=1,2,3,\dots $,
$$
A\biggl( \sum _{j=1}^m R_j\beta _j\biggr)  \left\{
\begin{aligned}
&\geq 0
& \quad \text{on } M\setminus W_m\\
&>1 &\quad \text {on } \overline W_0.
\end{aligned}
\right.
$$
Let $R_1\geq 1$. Given $R_1,\dots ,R_{m-1}>0$ with the above
property, using the fact that $A\beta _m>1$ on $\overline
W_{m-1}$, we get, for $R_m\gg 0$,
$$ A\biggl( \sum _{j=1}^m R_j\beta _j\biggr)
>1\qquad \text {on } \overline W_{m-1}.
$$
On $M\setminus (W_{m-1}\cup W_m)$ we have $A\beta _m\geq 0$ and
hence
$$
A\biggl( \sum _{j=1}^m R_j\beta _j\biggr)\geq A\biggl( \sum
_{j=1}^{m-1} R_j\beta _j\biggr)\geq 0
$$
On $\overline W_0$, the above middle expression, and hence
the expression on the left, is greater
than~$1$. Proceeding, we get the sequence $\seq Rm$. The sum $\sum
R_m\beta _m$ is locally finite in $X$ and the sequence of sets
$\seq Wm$ is locally finite in $X$, so the sum converges to a \fn
$\alpha $ with the required properties.
\end{pf*}

A slight modification of the proof of Theorem~1.1 gives the
following more general version:

\begin{thm}
Suppose $K$ is a \cpt subset of $M$ whose complement $M\setminus
K$ has no \rel \cpt \compsns , $\rho $ is a positive \cont \fn on
$M$, and $W$ is a \nbd of $K$ in $M$.  Then there exists a $\cinf
$ \fn $\vphi $ on $M$ \st
\begin{enumerate}
\item[(i)] $\vphi \geq 0$ and $A\vphi \geq 0$ on $M$,
\item[(ii)] $\vphi >\rho $ and $A\vphi >\rho$ on $M\setminus W$,
\item[(iii)] $\vphi \equiv 0$ on $K$,
\item[(iv)] $\vphi >0$ on $M\setminus K$, and
\item[(v)] For each \cpt set $E\subset M\setminus K$, we have
$A\vphi >\delta $ on $E$ for some constant $\delta =\delta (\vphi
,E)>0$.
\end{enumerate}
\end{thm}
\begin{pf} We proceed as in the proof of Theorem~1.1 but now with
$K_0=K$. By Lemma~1.5, we may choose nonempty \cpt sets $\seq K\nu
$ \st $M=\bigcup _{\nu =1}^{\infty }K_\nu $ and \stns , for
each~$\nu =0,1,2,\dots $, we have $K_\nu \subset \overset \circ
K_{\nu +1}$ and $M\setminus K_\nu $ has no \rel \cpt \compsns .

Given a point $p\in M\setminus K$, there is a unique $\nu =\nu
(p)\geq 0$ with $p\in K_{\nu +1}\setminus K_\nu $ and, by
Proposition~1.3, there is a $\cinf $ nonnegative \fn $\alpha _p$
and a \rel \cpt \nbd $V_p$ of~$p$ in $M\setminus K_\nu $ \st
$A\alpha _p\geq 0$ on $M$, $\alpha _p\equiv 0$ on $K _\nu $, and
$\alpha _p>\rho $ and $A\alpha _p>\rho $ on $V_p$. Thus we may
choose a sequence of points $\seq pk$ in $M\setminus K$ and
corresponding \fns $\seq \alpha {p_k}$ and \nbds $\seq V{p_k}$ so
that $\seq V{p_k}$ forms a covering of~$M\setminus W$ which is
locally finite in $M$ (as in the proof of Theorem~1.1, one may
take $\seq pk$ to be an enumeration of $\cup _{\nu =0}^\infty
Z_\nu $ where, for each $\nu $, $Z_\nu $ is a finite set of points
in $M\setminus  (W\cup \overset\circ K_\nu )$ \st $\{ V_p\} _{p\in
Z_\nu }$ covers $K_{\nu +1}\setminus (W\cup \overset \circ K_\nu )
$). The collection $\{ \supp \alpha _{p_k} \} $ is then locally
finite in $M$ and the locally finite sum $\sum _{k=1}^\infty
\alpha _{p_k}$ converges to a $\cinf $ \fn $\psi $ on $M$
satisfying $\psi \geq \alpha _{p_k}$ and $A\psi \geq A\alpha
_{p_k}$ for each $k$. It follows that $\psi \geq 0$ and $A\psi
\geq 0$ on $M$, $\psi \equiv 0$ on $K=K_0$, and $\psi >\rho $ and
$A\psi >\rho $ on $M\setminus W$.

In order to obtain the properties~(iv)~and~(v), we choose a
sequence of points $\seq qm$ in $M\setminus K$ and corresponding
\fns $\seq \alpha {q_m}$ and \nbds $\seq V{q_m}$ so that $\seq
V{q_m}$ covers $W\setminus K$. Applying a diagonal argument,
we may choose a sequence of positive
numbers $\seq \epsilon m $ converging to~$0$ so fast that each
derivative of arbitrary order for the sequence of partial sums
of $\sum \epsilon _m\alpha _{q_m}$ converges uniformly on \cpt
subsets of $M$. The \fn $\vphi=\psi +\sum \epsilon _m\alpha
_{q_m}$ will then have the required properties.
\end{pf}

The main topological fact required in the proof of Proposition~1.2
was part~(ii) of Lemma~1.9. This fact is slightly easier to verify
for $U$ a \comp of the complement of a \cpt set. But the more
general version (as stated in Lemma~1.9) and the proof of
Proposition~1.2 actually yield the following:

\begin{prop}
Let $U$ be a \con open subset of $M$ which contains a \con non\cpt
closed subset of $M$. Then, for each point $p\in U$, there exists
a $\cinf $ \fn $\alpha $ \st
\begin{enumerate}
\item[(i)] $\alpha \geq 0$ and $A\alpha \geq 0$ on $M$,
\item[(ii)] $\supp \alpha \subset U$,
\item[(iii)] $\alpha (p)>0$, and
\item[(iv)]  $A\alpha >1$ on a \nbd of $p$.
\end{enumerate}
\end{prop}

\begin{prop}
Let $K$ be a closed subset of $M$ \st each \comp of $M\setminus K$
contains a \con non\cpt closed subset of $M$. Then, for each point
$p\in M\setminus K$, there is a $\cinf $ nonnegative \fn $\alpha $
on $M$ \st $A\alpha \geq 0$ on $M$, $\alpha \equiv 0$ on $K$,
$\alpha (p)>0$, and $A\alpha
>1$ near~$p$.
\end{prop}

We also get a corresponding generalization of Theorem~1.10:

\begin{thm}
Suppose $K$ is a closed subset of $M$ \st each \comp of
$M\setminus K$ contains a \con non\cpt closed subset of $M$ and
$D\subset M\setminus K$ is a closed subset of $M$ with no \cpt
\compsns . Then, for every positive \cont real-valued \fn $\rho $
on $M$, there is a $\cinf $ \fnns~$\vphi $ \st
\begin{enumerate}
\item[(i)] $\vphi \geq 0$ and $A\vphi \geq 0$ on $M$,
\item[(ii)] $\vphi >\rho $ and $A\vphi >\rho$ on $D$,
\item[(iii)] $\vphi \equiv 0$ on $K$,
\item[(iv)] $\vphi >0$ on $M\setminus K$, and
\item[(v)] For each \cpt set $E\subset M\setminus K$, we have
$A\vphi >\delta $ on $E$ for some constant $\delta =\delta (\vphi
,E)>0$.
\end{enumerate}
\end{thm}

Before addressing the proof, we consider some remarks.

\begin{rmks}
\noindent 1. If the coefficients of $A$ are (for example) $C^1$
and the constant term is nonpositive, then the existence of a \con
non\cpt closed subset~$C$ of $M$ with $C\subset U$ is necessary in
Proposition~1.11. In fact, if $U$ is a \con open subset of $M$ and
$M$ admits a nonconstant nonnegative upper semi-\cont subsolution
$\alpha $ which vanishes on~$M\setminus U$, then $U$ must contain
such a set. For $\alpha (p)>0$ at some point $p\in U$ and hence we
may choose a number $\epsilon $ with $0<\epsilon <\alpha (p)$ and
a \nbd $V$ of the closed set $\setof{x\in M}{\alpha (x)\geq
\epsilon }$ with $\overline V\subset U$.
The maximum principle then
implies that the \comp $W$ of $V$ containing $p$ is not \rel \cpt
in $M$. Thus the set $C=\overline W$ is a closed \con non\cpt
subset of $M$ contained in $U$.

In particular, as the following example illustrates, the
conclusions of Proposition~1.2 and Proposition~1.3 do {\it not}
hold in general for $K$ a closed non\cpt set.
\end{rmks}

\begin{exmp}
Let $K$ be the closed subset of the manifold
$M=\R ^2\setminus \{ (0,0)\} $ given by
$$
K=\bigl( M \setminus (0,2)\times (0,2) \bigr) \cup \bigcup
_{m=1}^\infty \{ 1/m \} \times [0,1].$$ Then the complement
$U=M\setminus K$ is \con and $\overline U$ is non\cptns . But $U$
does not contain a \con non\cpt closed subset of $M$ and,
therefore, every nonnegative upper semi-\cont sub\harm \fn $\vphi
$ on $M$ which vanishes on $K$ must vanish everywhere in~$M$.
\end{exmp}

2. As the following example shows, the conclusion of Theorem~1.13
may fail to hold if the set $D\subset M\setminus K$ has \cpt
\compsns .

\begin{exmp}
The complement $U=M\setminus K$ in $M=\R ^2\setminus \{ (0,1)\} $
of the closed set
$$
K=\bigl( M\setminus (0,\infty )\times (0,4)\bigr) \cup \bigcup
_{m=1}^\infty \{ 1/(2m) \} \times [0,2]$$ is a \con set with
non\cpt closure and $U$ contains the closed non\cpt \con set
$C=[1,\infty )\times \{ 3 \} $. The non\cpt subset
$$
D=\setof {(1/(2m+1),1)}{m\in \N }$$ of $U$ is closed (in fact,
discrete) in $M$. If $\vphi $ is a nonnegative upper semi-\cont
 sub\harm \fn on $M$ which vanishes on~$K$, then, for each~$m$,
applying the maximum principle in $[1/(2m+2),1/(2m)]\times [0,2]$,
we get a number $\seq rm$ with $1/(2m+2)<r_m<1/(2m)$ and $\vphi
(r_m,2)\geq \vphi (1/(2m+1),1)$. Since $(r_m,2)\to (0,2)\in
K\subset M$ and $\vphi $ is upper semi-continuous, it follows that
$\vphi $ must be bounded on~$D$.
\end{exmp}

\noindent 3. If $K\subset M$ is a \cpt set, then one can achieve
the conditions in Proposition~1.3 and Theorem~1.10 by replacing
$K$ by the \cpt set~$\widehat K$. Because we have Proposition~1.12
and Theorem~1.13, for a general closed set $K\subset M$ it is
natural to define $\widehat K$ to be the union of $K$ with all
\comps of $M\setminus K$ which do {\it not} contain any \con
non\cpt closed subsets of $M$.

The main step in the proof of Theorem~1.13 is the case in which
$D$ and $M\setminus K$ are \conns .

\begin{lem}
Suppose $U$ is a \con open subset of $M$, $C$ is a \con non\cpt
closed subset of $M$ with $C\subset U$, and $\rho $ is a positive
\cont \fn on $M$. Then there is a $\cinf $ \fnns~$\vphi $ \st
\begin{enumerate}
\item[(i)] $\vphi\geq 0$ and $A\vphi \geq 0$ on $M$,
\item[(ii)] $\vphi >\rho $ and $A\vphi >\rho $ on $C$,
\item[(iii)] $\vphi\equiv 0$ on $M\setminus U$,
\item[(iv)] $\vphi >0$ on $U$, and
\item[(v)] For each \cpt set $E\subset U$,
we have $A\vphi >\delta $ on $E$ for some $\delta =\delta (\vphi , E)>0$.
\end{enumerate}
\end{lem}

\begin{pf}
We first show that there is a nonnegative $\cinf $ \fn $\psi $ \st
$A\psi \geq 0$ on $M$, $\psi\equiv 0$ on $M\setminus U$, and $\psi
>\rho $ and $A\psi >\rho $ on~$C$.

For this purpose, we may assume without loss of generality that
$C$ is locally \conns . For we may choose (as in the proof of
Lemma~1.9) a locally finite (in~$M$) covering $\mathcal A$ of $C$
by \rel \cpt \con open subsets of $U$.  We may also choose the
covering so that each element meets $C$ and has locally connected
closure (for example, we may choose~$\mathcal A$ so that, for each
$B\in \mathcal A$, there is a diffeomorphism of some \nbd of
$\overline B$ onto an open subset of $\R ^n$ mapping $B$ onto a
ball). The closed \con non\cpt set
$$
C'\equiv \overline {\bigcup _{B\in \mathcal A}B} =\bigcup _{B\in
\mathcal A}\overline B \subset U
$$
is then locally \conns . For if $p\in C'$ and $B_1,\dots , B_k$
are the (finitely many) elements of $\mathcal A$ whose closures
contain $p$, then, for each $j=1,\dots ,k$, we may choose a \nbd
$W_j$ of~$p$ in $M$ \st $W_j\cap \overline B_j$ is \con and
$W_j\cap\overline B=\emptyset $ for each set $B\in \mathcal
A\setminus \{ \, B_1,\dots , B_k\, \} $. The set $D\equiv \bigcup
_{j=1}^k(W_j\cap \overline B_j)$ is then a \con subset of $C'$
which contains the set $W_1\cap \cdots \cap W_k\cap C'$, a \nbd of
$p$ relative to~$C'$.  It follows that $C'$ is locally \con
(since, by choosing the \nbds $\seq Wj$ small, one sees that the
\comps of any open subset of $C'$ are open relative to $C'$).
Therefore, by replacing $C$ with the set $C'$, we may assume that
$C$ is locally \conns .

By Lemma~1.4, there is a sequence of \cpt sets $\seq K\nu $ \st
$M=\bigcup _{\nu =1}^\infty K_\nu $ and, for each~$\nu $, we have
$K_\nu \subset \overset \circ K_{\nu +1}$ and $C\setminus K_\nu $
has only finitely many \compsns , all of which have non\cpt
closure. For we may choose inductively a sequence of \cpt subsets
$\{ K'_\nu \} $ of $M$ \st $M=\bigcup _{\nu =1}^\infty K'_\nu $
and \stns , for each~$\nu $, we have $$ K_\nu \equiv K'_\nu \cup
\widehat {(K'_\nu \cap C)}_C\subset\overset {\circ \,\,}{K'}_{\nu
+1};$$ where, for $K\subset C$ \cptns , $\widehat K_C$ is the
union of $K$ with all of the \rel \cpt \comps of $C\setminus K$.

We now proceed as in the proofs of Theorem~1.1 and Theorem~1.10.
Let $K_0=\emptyset $. Given $p\in C$, there is a
unique $\nu =\nu (p)$ with $p\in
K_{\nu +1}\setminus K_\nu $. The \comp $U_p$ of $U\setminus K_\nu
$ containing~$p$ must also contain the closure of some \comp of
$C\setminus K_{\nu +1}\subset C\setminus \overset \circ K_{\nu +1}
\subset C\setminus K_\nu $. For we may take a point $q$ in the
\comp of $C\setminus K_\nu $ containing~$p$ (a set with non\cpt
closure) which lies outside~$K_{\nu +1}$. The closure of the \comp
of $C\setminus K_{\nu +1}$ containing~$q$ is then contained
in~$U_p$.  We may apply Proposition~1.11 to get a $\cinf $
nonnegative \fn $\alpha _p$ and a \rel \cpt \nbd $V_p$ of $p$ in
$U_p$ \st $A\alpha _p\geq 0$ on~$M$, $\supp \alpha _p\subset U_p$,
and $\alpha _p>\rho $ and $A\alpha _p>\rho $ on $V_p$. Thus we may
choose a sequence of points $\seq pk$ in $C$ and corresponding
\fns $\seq \alpha {p_k}$ and \nbds $\seq V{p_k}$ so that $\seq
V{p_k}$ forms a locally finite (in $M$) covering of~$C$. The
collection $\{ \supp \alpha _{p_k} \} $ is then locally finite in
$M$ because $\supp \alpha _{p_k}\subset U\setminus K_\nu $
whenever $p_k\not\in K_\nu $. Hence the sum $\sum _{k=1}^\infty
\alpha _{p_k}$ is locally finite and, therefore, convergent to a
$\cinf $ \fn $\psi $ on $M$ with the required properties.

Now, by Lemma~1.9, there is a sequence of \con open sets $\seq
U\nu $ \st $U=\bigcup _{\nu =1}^\infty U_\nu $, $C\subset U_1$,
and, for each~$\nu $, $\overline U_\nu \subset U_{\nu +1}$.  By
the above, we may form a $\cinf $ nonnegative \fn $\psi $ \st
$A\psi \geq 0$ on $M$, $\psi \equiv 0$ on $M\setminus U_1$, and
$\psi
>\rho $ and $A\psi  >\rho $ on $C$ and, for each
$\nu =1,2,3,\dots $, we may form a $\cinf $ nonnegative \fn $\psi
_\nu $ \st $A\psi _\nu \geq 0$ on $M$, $\psi _\nu \equiv 0$ on
$M\setminus U_{\nu +1}$, and $\psi _\nu
>1$ and $A\psi _\nu >1$ on $\overline U_\nu $. Choosing a
sequence of positive numbers $\seq \epsilon \nu $ converging
to~$0$ sufficiently fast, the \fn $\vphi=\psi +\sum \epsilon _\nu
\psi _\nu $ will have the required properties.
\end{pf}

For the general case, we will apply the following:
\begin{lem}
Suppose $X$ is a second countable, \conns , locally \conns ,
locally \cptns , Hausdorff topological space; $K$ is a closed
subset of~$X$; and $D\subset X\setminus K$ is a closed subset of
$X$ with no \cpt \compsns . Then there exists a countable locally
finite (in $X$) family of disjoint \con non\cpt closed sets $\seq
C\lambda _{\lambda \in \Lambda }$ and a locally finite (in $X$)
family of disjoint
\con open sets $\seq U\lambda
_{\lambda \in \Lambda }$ \st
$$
D\subset C\equiv \bigcup _{\lambda \in \Lambda }C_\lambda
\qquad\text{and}\qquad C_\lambda \subset U_\lambda \subset
\overline U_\lambda \subset X\setminus K \quad \forall\, \lambda
\in \Lambda .
$$
\end{lem}
\begin{rmk}
We will not use the fact that the sets $\seq U\lambda $ are
disjoint.
\end{rmk}
\begin{pf}
As in the proof of Lemma~1.9, there is a countable locally finite
covering $\mathcal A_D$ of $D$ by \con open \rel \cpt subsets of
$X\setminus K$  which meet $D$. Thus
$$
D\subset V\equiv \bigcup _{B\in \mathcal A_D}B\subset \overline
V=\bigcup _{B\in \mathcal A_D}\overline B\subset X\setminus K
$$
(where we have used the local finiteness of the collection
$\mathcal A_D$). Since each of the \comps of $V$ meets, and
therefore contains, a \comp of~$D$, the family of \comps $\seq
V\gamma _{\gamma \in \Gamma }$ of $V$ is a locally finite family
of \con open sets with non\cpt closure. The set
$$
C\equiv \overline V=\bigcup _{\gamma \in \Gamma }\overline
V_\gamma
$$
is a closed set contained in $X\setminus K$ and the family of \comps
$\seq C\lambda _{\lambda \in \Lambda }$ of $C$ satisfies
$$
C_\lambda =\bigcup_ {\gamma \in \Gamma, \, \overline V_\gamma
\subset C_\lambda }\overline V_\gamma =\overline{\bigcup_ {\gamma
\in \Gamma, \, \overline V_\gamma \subset C_\lambda }V_\gamma }
\quad \forall \, \lambda \in \Lambda .
$$
It follows that the family is locally finite in $X$ (since the
family $\seq {\overline V}\gamma _{\gamma \in \Gamma }$ is locally
finite) and that $C_\lambda $ is closed for each $\lambda \in
\Lambda $. Consequently, we may choose a locally finite covering
$\mathcal A_C$ of $C$ by \con \rel \cpt open subsets of
$X\setminus K$ \stns , for
each element $B\in \mathcal A_C$,
$B$ and $\overline B$ meet exactly one \comp of $C$.
For each $\lambda \in \Lambda $, taking $U_\lambda $ to be
the \comp of the set
$$
\bigcup _{B\in \mathcal A_C, B\cap C_\lambda \neq
\emptyset }B \setminus
\bigcup _{B\in \mathcal A_C, B\cap C_\lambda =
\emptyset }\overline B
$$
containing $C_\lambda $,
we get disjoint
\con open sets $\seq U\lambda _{\lambda \in
\Lambda }$ with
$C_\lambda \subset U_\lambda \subset \overline U_\lambda
\subset X\setminus K$
for each $\lambda \in \Lambda $. This family is locally finite in
$X$. For each point in $X$ has a \nbdns~$Q$ which meets only
finitely many elements $B_1,\dots , B_k$ of $\mathcal A_C$. Each
$B_j$ meets a unique \comp $C_{\lambda _j}$ of $C$. If $\lambda
\in \Lambda $ with $Q\cap U_\lambda \neq \emptyset $, then $Q\cap
B\neq \emptyset $ for some $B\in \mathcal A_C$ with $B\cap
C_\lambda \neq \emptyset $. Hence we must have $B=B_j$ for
some~$j$ and, therefore, $\lambda =\lambda _j$.
\end{pf}

\begin{pf*}{Proof of Theorem~1.13}
Let
$$ D\subset C\equiv \bigcup _{\lambda \in \Lambda }C_\lambda
\qquad\text{and}\qquad C_\lambda \subset U_\lambda \subset
\overline U_\lambda \subset M\setminus K \quad \forall\, \lambda
\in \Lambda
$$
be as in Lemma~1.17. Applying Lemma~1.16 to each pair of sets
$C_\lambda \subset U_\lambda $, we get a nonnegative $\cinf $ \fn
$\alpha _\lambda $ \st $A\alpha _\lambda \geq 0$ on $M$, $\alpha
_\lambda \equiv 0$ on $M\setminus U_\lambda $, and $\alpha
_\lambda >\rho $ and $A\alpha _\lambda >\rho $ on $C_\lambda $ (we
do not need the conditions~(iv)~and~(v) of Lemma~1.16 for this
part). Since the family $\seq U\lambda $ is locally finite in $M$,
the sum $\sum \alpha _\lambda $ determines a nonnegative $\cinf $
\fn $\alpha $ with $A\alpha \geq 0$ on $M$, $\alpha \equiv 0$ on
$M\setminus \bigcup _{\lambda \in \Lambda }U_\lambda \supset K$,
and $\alpha >\rho $ and $A\alpha >\rho $ on $C\supset D$.

Applying Lemma~1.16 to each of the \comps ${\seq Vj}_{j\in J}$ of
$M\setminus K$, we get, for each $j\in J$, a $\cinf $ nonnegative
\fn $\beta _j$ \st $A\beta _j\geq 0$ on $M$, $\beta _j \equiv 0$
on $M\setminus V_j$, $\beta _j>0$ on $V_j$, and, in $V_j$, $A\beta
_j$ is locally bounded below by a positive constant. For $J$ a
finite set, we may now take $\vphi =\alpha +\sum _{j\in J}\beta
_j$. If $J$ is infinite, then, assuming as we may that $J=\N $ and
choosing a sequence of positive numbers $\seq \epsilon j$
converging to~$0$ sufficiently fast, the \fn $\vphi =\alpha +\sum
_{j=1}^\infty \epsilon _j\beta _j$ will have the required
properties.
\end{pf*}

We close this section with the following observation concerning
Theorem~1.10 for $K$ the closure of a smooth \rel \cpt domain.
\begin{cor}
Suppose $\Omega $ is a $\cinf $ \rel \cpt domain in~$M$ whose
complement $M\setminus \Omega $ has no \cpt \compsns , $\rho $ is
a positive \cont \fn on $M$, and $W$ is a \nbd of $\overline
\Omega $. Then there is a $\cinf $ \fnns~$\vphi $ on $M$ \st
$\vphi >\rho $ on $M\setminus W$, $A\vphi >\rho $ on $M$, $0$ is a
regular value for~$\vphi $, and
$\Omega =\setof {x\in M}{\vphi (x)<0}$.
\end{cor}
\begin{pf}
There exists a $\cinf $ \fn $\tau $ on $M$ \st $0$ is a regular
value for~$\tau $, $\tau $ is locally constant on $M\setminus V$
for some \rel \cpt \nbd $V$ of $\partial \Omega $ in $W$, and $
\Omega =\setof {x\in M}{\tau (x)<0}$. For $\epsilon >0$
sufficiently small, we have
$$
D\equiv\setof{x\in M}{-2\epsilon \leq \tau (x)\leq 2\epsilon }
\subset V
$$
and $(d\tau )_x\neq 0$ for each point $x\in D$. By Theorem~1.1,
there is a $\cinf $ \fn $\alpha $ with \cpt support in $\Omega $
\st $\alpha \leq 0$ on $M$ and $A\alpha >\rho $ on $\setof{x\in
M}{ \tau (x)\leq -\epsilon }$. By Theorem~1.10, there is a $\cinf
$ nonnegative \fnns~$\beta $ on~$M$ \st $A\beta \geq 0$ on $M$,
$\beta \equiv 0$ on
$\Omega $, $\beta >0$ and $A\beta >0$ on $M\setminus \overline
\Omega $, and $\beta >\rho $ and $A\beta >1+\rho $ on $\setof{x\in
M}{\tau (x)\geq \epsilon }$.
Finally, we may fix a $\cinf $ \fn $\chi :\R \to
[0,\infty )$ as in Lemma~1.8. Let $R_1,R_2,R_3>1$ and let
$$
\vphi =\alpha +R_2\chi \bigl( R_1(\tau +2\epsilon )\bigr) -R_2\chi
(2R_1\epsilon )+R_3\beta .
$$
On $M\setminus W$, we have $\vphi \geq R_3\beta >\beta >\rho $.
On $\Omega $, we have
$$
\vphi <0+R_2\chi \bigl( R_1(0 +2\epsilon )\bigr) -R_2\chi
(2R_1\epsilon )+R_3\cdot 0 =0.
$$
On $M\setminus \overline \Omega $, we have
$$
\vphi \geq 0+R_2\chi \bigl( R_1(0 +2\epsilon )\bigr) -R_2\chi
(2R_1\epsilon )+R_3\beta \geq R_3\beta >0.
$$
Thus $\Omega =\setof {x\in M}{\vphi (x)<0}$. For any point $x\in
\partial \Omega =\setof {x\in M}{\vphi (x)=0}$,
we have $\alpha =0$ near~$x$ and $\beta $ has a local minimum
at~$x$. Thus
$$
d\vphi =d\alpha +R_1R_2\chi '(2R_1\epsilon )d\tau +R_3d\beta
=0+R_1R_2\chi '(2R_1\epsilon )d\tau +0\neq 0.
$$
By Lemma~1.7, for $R_1\gg 0$, we get $A\bigl[ \chi \bigl( R_1(\tau
+2\epsilon )\bigr) \bigr] \geq 0$ on $\Omega \cup D=\setof {x\in
M}{\tau (x)\leq 2\epsilon }$ and $A\bigl[ \chi \bigl( R_1(\tau
+2\epsilon )\bigr) \bigr]$ is bounded below by a positive constant
in a \nbd of each point in $\setof {x\in M}{-2\epsilon <\tau
(x)\leq 2\epsilon }\subset D$. On $\setof {x\in M}{\tau (x)\leq
-\epsilon }$ we have
$$
A\vphi \geq A\alpha +0+0=A\alpha >\rho .
$$
For $R_2\gg 0$, on $\setof {x\in M}{-\epsilon \leq \tau (x)\leq
\epsilon }$ we have
$$
A\vphi \geq A\alpha +R_2A\bigl[ \chi \bigl( R_1(\tau +2\epsilon
)\bigr) \bigr] > \rho .
$$
Finally, since $\tau $ is locally constant on $M\setminus V$, for
$R_3\gg 0$, on $\setof {x\in M}{\tau (x)\geq \epsilon }\subset
M\setminus \overline \Omega $ we have
$$
A\vphi =0 +R_2A\bigl[ \chi \bigl( R_1(\tau +2\epsilon )\bigr)
]+R_3A\beta > R_2A\bigl[ \chi \bigl( R_1(\tau +2\epsilon )\bigr)
]+R_3(1+\rho )>\rho .
$$
\end{pf}

\section{Two applications}

To illustrate the broad utility of the existence of exhausting
strict subsolutions,  we consider two (mostly known) consequences.

We first recall that any $C^k$ \fn $\vphi $ on a smooth manifold
can be approximated in the $C^k$ Whitney topology by a $\cinf $
Morse \fnns~$\psi $ [GG].  Applying this to a \fn $\vphi $ from
Theorem~1.1, we get
\begin{cor}
If $M$ and $A$ are as in Section~1 with $M$ non\cpt and $\rho $ is
a positive \cont \fn on $M$, then there exists a $\cinf $ Morse
\fn $\psi $ satisfying $\psi >\rho $ and $A\psi
>\rho$.
\end{cor}

Taking $A$ to be $\lap _g$ for a Riemannian metric $g$ and $\rho $
to be a \cont \exh \fnns , we get a $\cinf $ \exh \fn $\psi $
whose Hessian is of index at most~$n-1$
at each point. Thus we get the following well-known fact:
\begin{thm}
A \con non\cpt $\cinf $ manifold of dimension~$n$ has the homotopy
type of a CW~complex with cells of dimension~$\leq n-1$.
\end{thm}

The next observation is that the existence of exhausting strict
subsolutions allows one to construct a Hermitian metric of
positive scalar curvature (positive curvature in the case of a
Riemann surface) in a \holo line bundle with a nontrivial \holo
section.  In the \cpt case, such Hermitian
\holo line bundles are, of course, a natural substitute for \str
sub\harm \fnsns .

For the rest of this section, $X$ will denote a \con \cpx manifold
of (\cpxns ) dimension~$n$ and $g$ will denote a $\cinf $
Hermitian metric in~$X$. The {\it Levi form} of a $C^2$
\fnns~$\vphi $ on~$X$ is the Hermitian tensor given by, in local
\holo coordinates $(z_1,\dots ,z_n)$,
$$
\lev {\vphi } =\sum _{i,j=1}^n \frac {\partial ^2\vphi }{\partial
z_i\partial \bar z_j} dz_i d\bar z_j.
$$
The {\it Laplace operator} $\lap _g$ for the Hermitian metric~$g$
is given by the trace of the Levi form:
$$
\lap _g=\sum _{i,j} g^{i\bar j}\frac {\partial ^2}{\partial
z_i\partial \bar z_j}
$$
where $\bigl( g^{i\bar j} \bigr) =\overline {\bigl( g_{i\bar j}
\bigr)\inv }$. This elliptic operator is equal to $1/2$ the
Laplace operator of the associated Riemannian metric if $g$ is
K\"ahler. A $C^2$ real-valued \fnns~$\vphi $ is called {\it
subharmonic} ({\it \str subharmonic}) with respect to~$g$
if $\lap _g\vphi \geq 0$
(respectively, $\lap _g\vphi > 0$).

If $L$ is a \holo line bundle on $X$ and $h$ is a $C^2$ Hermitian
metric in~$L$, then the {\it curvature} of $h$ is the Hermitian
tensor $\Theta _h$ given by
$$
\Theta _h\equiv \lev {-\log |s|^2_h}
$$
for any nonvanishing local \holo section $s$ of $L$.  The {\it
scalar curvature}~$\mathcal R_h$ of $h$ with respect to~$g$ is
given by the trace of the curvature; that is, locally,
$$
\mathcal R_h\equiv \lap _g (-\log |s|^2_h).
$$
In particular, if $X$ is a Riemann surface, then $\mathcal
R_h=\Theta _h/g$.

\begin{thm}
Let $L$ be a \holo line bundle on~$X$.  If $X$ is non\cpt or
$L=[D]$ is the \holo line bundle associated to a nontrivial
effective divisor~$D$ in~$X$ (i.e.~$L$ is a \holo line bundle
which admits a nontrivial global \holo section), then $L$ admits a
$\cinf $ Hermitian metric~$h$ with positive scalar curvature.
\end{thm}
\begin{pf}
Fix a $\cinf $ Hermitian metric $k$ in $L$. We will modify $k$ to
obtain~$h$.

Assuming first that $X$ is non\cptns , Theorem~1.1 provides a
$\cinf $ \str sub\harm (with respect to~$g$) \exh \fnns~$\vphi $.
If $\chi $ is a $\cinf $ \fn on~$\R $ with $\chi '>0$ and $\chi
''\geq 0$ and
$$
h=e^{-\chi (\vphi )}k,
$$
then
$$
\mathcal R_h=\lap _g(\chi (\vphi ))+\mathcal R_k \\
=\chi '(\vphi )\lap _g\vphi + \chi ''(\vphi )|\partial \vphi
|^2_g +\mathcal R_k\\
\geq \chi '(\vphi )\lap _g\vphi +\mathcal R_k.
$$
Choosing $\chi $ so that $\chi '(t)\to \infty $ sufficiently fast
as $t\to \infty $, we get $\mathcal R_h>0$.

Assuming now that $X$ is \cpt and $L=[D]$, where $D$ is a
nontrivial effective divisor, let $Y=|D|\subset X$ be the support
of~$D$ and let~$s$ be a global \holo section of~$L$ with
associated divisor~$D$.  Applying Theorem~1.1 to a non\cpt \nbd
of~$Y$ in~$X$ and cutting off, we get a $\cinf $ \fnns~$\alpha $
on~$X$ which is \str sub\harm on a \nbdns~$U$ of $Y$.  After
shrinking $U$ slightly and replacing $\alpha $ by a large
multiple, we may assume that
$$
\lap _g\alpha +\mathcal R_k>2 \qquad \text{on } U.
$$
Since $-\log |s|^2_k\to \infty $ at~$Y$, we have, for $N\gg 0$,
$$
Y\subset \setof {x\in X}{\alpha (x)-\log |s(x)|^2_k\geq N}\subset
U
$$
(setting $\alpha -\log |s|^2_k=\infty $ along $Y$).  We may choose
a $\cinf $ \fn $\lambda $ on~$\R $ \st $\lambda '\geq 0$, $\lambda
''\geq 0$, $\lambda (t)=t$ if $t\geq 3N$, and $\lambda (t)=2N$ if
$t\leq N$. We set $\lambda (\infty )=\infty $. The restriction of
the \fn
$$
\rho \equiv \lambda \bigl(\alpha -\log |s|^2_k\bigr)
$$
to $X\setminus Y$ is $\cinf $ and sub\harm because $\rho \equiv
2N$ on a \nbd of $X\setminus U$ (and hence $\lap _g\rho =0$),
while on $U\setminus Y$ we have
$$
\lap _g\rho = \lambda '\bigl(\alpha -\log |s|^2_k\bigr) \cdot
\bigl(\lap _g\alpha +\mathcal R_k\bigr) +\lambda ''\bigl(\alpha
-\log |s|^2_k\bigr)\cdot \bigl|\partial \bigl(\alpha -\log
|s|^2_k\bigr)\bigr| ^2_g \geq 2\lambda '\bigl(\alpha -\log
|s|^2_k\bigr)\geq 0.
$$
Observe also that $\rho =\alpha -\log |s|^2_k$ on the \rel \cpt
\nbdns~$V$ of $Y$ in $U$ given by
$$
V=\setof{x\in X}{\alpha (x)-\log |s(x)|^2_k>3N}.
$$
Applying Theorem~1.1 to the \con non\cpt manifold $X\setminus Y$
and cutting off near $Y$, we get a $\cinf $ \fnns~$\beta $ with
\cpt support in $X\setminus Y$ satisfying $\lap _g\beta >0$ on
$X\setminus V$. Choosing $\epsilon >0$ so small that $\epsilon
\lap _g\beta
>-1$ on $X$, we see that the restriction of the \fn $\gamma
\equiv \rho +\epsilon \beta $ to $X\setminus Y$ is $\cinf $ and
\str sub\harmns . In fact, on $V\setminus Y$, we have
$$
\lap _g\gamma =\lap (\alpha -\log |s|^2_k +\epsilon \beta )
>2-1=1.
$$
We may now define $|\xi |_h^2$ for $\xi \in
L_x$ with $x\in X$ by
$$
|\xi |_h^2 = \left\{
\begin{aligned}
&e^{-\gamma (x)}|\xi /s(x)|^2
& \quad \text{if } x\in X\setminus Y\\
&e^{-\alpha (x)-\epsilon\beta (x)}|\xi |_k^2&\quad \text {if }
x\in V
\end{aligned}
\right.
$$
Then $h$ is a well-defined $\cinf $ Hermitian metric in~$L$ since,
for $x\in V\setminus Y$ and $\xi \in L_x$, we have
$$
e^{-\gamma (x)}|\xi /s(x)|^2=e^{-\bigl(\alpha (x)-\log
|s(x)|^2_k+\epsilon \beta (x)\bigr)}|\xi |_k^2/|s(x)|_k^2
=e^{-\alpha (x)-\epsilon \beta (x)}|\xi |_k^2.
$$
Furthermore, on $X\setminus Y$ we have
$$
\mathcal R_h=\lap _g \bigl( -\log |s|_h^2\bigr) =\lap _g\gamma
\left\{
\begin{aligned}
&>0
& \quad \text{on } X\setminus Y\\
&>1&\quad \text {on } V\setminus Y
\end{aligned}
\right. $$ By continuity, we also have $\mathcal R_h \geq 1>0$ at
points in $Y$. Thus $\mathcal R_h>0$ on $X$.
\end{pf}

For $X$ a Riemann surface, the above proofs become especially
simple.  For example, the construction of $\alpha $ in the proof
of Theorem~2.3 is trivial for $\dim X=1$ because $Y$ is discrete.
For $X$ an open Riemann surface, one gets a $\cinf $ \str \plsh
\exh \fn and, therefore, by [G] and [DG], one gets the theorem of
[BS] that an open Riemann surface is Stein. For a general Riemann
surface, Theorem~2.3 yields the following familiar fact:
\begin{cor}
The \holo line bundle associated to a nontrivial effective divisor
in a Riemann surface $X$ admits a $\cinf $ Hermitian metric with
positive curvature.
\end{cor}

\bibliographystyle{amsalpha.bst}

\end{document}